Oct. 15, 2017.

# On Some Sums of Digamma and Polygamma Functions – Revision (2017) and Review

Michael Milgram*,

Consulting Physicist, Geometrics Unlimited, Ltd.,

Box 1484, Deep River, Ont. Canada. K0J 1P0

**Abstract:** This paper is an enhanced version of a more than decade-older paper with a similar title. Many formulae involving both finite and infinite sums of digamma and polygamma functions up to quadratic order, few of which appear in standard reference works or the literature, but which periodically arise in applications, are collected, reviewed, listed and developed to the point that a knowledgeable reader could devise a formal proof. Several errors in the literature are corrected.

**Author's Note (2017 version):** The forerunner of this paper was originally consigned to arXiv in 2004 after being rejected by several reputable journals, mostly on the grounds that the contents were well-known, although one referee felt it necessary to opine that the author lacked a "prestigious affiliation". Over the intervening years, the arXiv entry (https://arxiv.org/abs/math/0406338) has garnered a reasonable number of citations in the published literature, thereby demonstrating a need for an easily accessible reference work, and the author has attempted to improve the prestige of his affiliation. Also, since that time, a number of results, particularly for finite sums, have been obtained that may be either not-so-well-known or new, many of which are needed for applications (e.g. [**31**]). At the same time, new results for Euler (harmonic) sums have independently appeared in the literature, but the correspondence between two equivalent results, one expressed as a sum of (powers of) harmonic numbers, the other as a sum of (powers or products of) digamma or polygamma functions is not always transparent, recognizable or easily determined. Thus, I have prepared this enhanced version of the original 2004 paper, containing a large number of additional, possibly new, results, intended to be used as a reference or review. As before, the emphasis has been on generalized Euler-type sums – that is, sums where the summand contains at least one digamma or polygamma function, and one or more independent parameters, but contains a minimal number (or no) Gamma-functions, except as needed for derivations. With specialized applications in mind, the results are limited to quadratic forms, although a few cubic forms are included when needed for simplification. Transcriptions of all the new additions have been checked by copying each formula from the text by hand and verifying numerically using the Maple computer code.

Note that the equation numbers in this and the original version do not always match, although the main text – Sections 1 to 3 – has only been modified slightly. In particular, Appendix B is now devoted to results for finite series. It is interesting to note that some results are closed, rather than consisting of transformations between series. A new Appendix C (mostly) repeats the previous listings of infinite series with some additions, again with only hints at the proofs. Appendix D consists of a computer-generated result that was too complicated to transcribe accurately and Appendix E lists a few simple specific and useful results extracted from the more general entries and the literature, regarding which, several corrections are noted. References to now-long-ago-published results, cited previously as unpublished, have also been updated. In general, within sections, results are listed in order of increasing complexity with only the barest of hints at proofs, except for some cases which are grouped together so that (implied) proofs might become more transparent.

* mike@geometrics-unlimited.com

**Abstract(2004):** Many formula involving sums of digamma and polygamma functions, few of which appear in standard reference works or the literature, but which periodically arise in applications, are collected and developed. Along the way, a new evaluation for some members of the family of hypergeometric functions – ${}_4F_3(1)$ – is presented, and a connection is made with Euler sums.

## 1. Introduction

Although sums of digamma functions often appear in applications such as particle transport [**19**] and the evaluation of Feynman diagrams [**7**,**24**,**25**] there is a paucity of results available in the standard references [**12**, **26**] or in the literature [**3**,**9**,**13**, **24**, **27**] except for the special case of Euler or harmonic sums [e.g. **10**]. Therefore, many formulae are often revisited by different authors [e.g. **27**] because no repository of quotable results exists. At the same time, many such sums are very slow to converge, so it is important to obtain closed form expressions when numerical evaluation is intended. The purpose of this paper is to collect and develop, in one place, a number of sums involving digamma and polygamma functions that have arisen in applications [e.g. **19** and **31**] or which were obtained as a byproduct of those endeavours. Although many of the results about to be quoted are easily established by evaluating parametric derivatives of hypergeometric functions of unit argument - a technique that "traces back to Newton" [**6**] - it is believed that many of them, which generalize previously known results, are also new. To keep the paper to a reasonable length, proofs are only sketched, sufficiently that a knowledgeable reader, with access to a computer algebra program, should be able to reproduce them.

By way of review, the digamma function $\psi(x)$ is defined by

$$\psi(x) = \frac{d}{dx}\log(\Gamma(x)) \qquad \textbf{(1)}$$

so that

$$\frac{d}{dx}\Gamma(x) = \psi(x)\Gamma(x) \qquad \textbf{(2)}$$

with*

$$\psi(x) = \sum_{l=1}^{k}\frac{1}{(x-l)} + \psi(x-k),\ k>0. \qquad \textbf{(3)}$$

Polygamma functions are written as:

$$\psi'(x) = \frac{d}{dx}\psi(x) \text{ and } \psi^{(n)}(x) \equiv \frac{d^n}{dx^n}\psi(x) \text{ if } n \ge 2 \qquad \textbf{(4)}$$

* see also (**7**).

Two useful, well-established results, used throughout, are

$$\sum_{\ell=0}^{k} \frac{1}{(\ell+a)} \frac{1}{(\ell+b)} = (\psi(k+a+1) - \psi(a) - \psi(k+b+1) + \psi(b))/(b-a) \quad \textbf{(5)}$$

and

$$\begin{aligned} \sum_{\ell=0}^{\infty} \frac{1}{(\ell+a)} \frac{1}{(\ell+b)} &= (\psi(a) - \psi(b))/(a-b) \qquad b \neq a \\ &= \psi'(a) \qquad if \qquad b = a \end{aligned} \quad \textbf{(6)}$$

Also, define generalized harmonic numbers [**27**, **6**, **35**]

$$H_n(a) \equiv \sum_{l=1}^{n} \frac{1}{(l+a-1)} = \psi(n+a) - \psi(a) \quad \textbf{(7)}$$

together with

$$H'_n(a) \equiv \sum_{l=1}^{n} \frac{1}{(l+a-1)^2} = -(\psi'(n+a) - \psi'(a)). \quad \textbf{(8)}$$

In general, if $k > 0$, commonly used variants are

$$H_n^{(k)}(1) \equiv \sum_{l=1}^{n} \frac{1}{l^k} = \frac{(-1)^{k+1}}{\Gamma(k)} \left(\psi^{(k-1)}(n+1) - \psi^{(k-1)}(1)\right), \quad \textbf{(9)}$$

$$H_n^{(k)}(\tfrac{1}{2}) \equiv \sum_{l=1}^{n} \frac{1}{(l-\frac{1}{2})^k} = \frac{(-1)^{k+1}}{\Gamma(k)} \left(\psi^{(k-1)}(n+\tfrac{1}{2}) - \psi^{(k-1)}(\tfrac{1}{2})\right) \quad \textbf{(10)}$$

and, with respect to the corresponding alternating series,

$$\begin{aligned} L_n^{(k)}(1) \equiv \sum_{l=1}^{n} \frac{(-1)^l}{l^k} = \frac{(-1)^k}{\Gamma(k)} \Big( &(-1)^n [\psi^{(k-1)}(n+1) - 2^{(1-k)} \psi^{(k-1)}(\tfrac{1}{2}n+1)] + (2^{(1-k)} - 1)\psi^{(k-1)}(1) \\ &+ 2^{(1-k)}(1-(-1)^n) \ln(2) \delta_{1,k} \Big) . \end{aligned} \quad \textbf{(11)}$$

Throughout, $\zeta(n)$ is Riemann's zeta function, $\delta_{1,k}$ is the Kronecker delta, $\zeta(s,c)$ is the Hurwitz zeta function, $\psi(1) = -\gamma$ is the Euler-Mascheroni constant and $Li_n(x) \equiv \sum_{l=1}^{\infty} \frac{x^l}{l^n}$ is the $n^{th}$ polylog function. Symbols starting with the letters "*j*" through "*n*" are always non-negative integers. All other symbols are continuous and (arbitrarily) complex. $\sum \equiv \sum_{l=0}^{\infty}$ and $\sum_{l=n_1}^{n_2} \cdots = 0$ if $n_2 < n_1$ (except see (**B.22**)). $G$ is Catalan's constant and $\lfloor x \rfloor$ symbolizes the "floor" function (greatest integer less than or equal to x).

## 2. Hypergeometric Sums

**Lemma 2.1.**

An easily obtained and useful sum needed later (see **(63)**), which generalizes a well-known result [**18**, Eq. 3.13(43), **26**, Eqs. 7.4.4.40-46], which is not listed in any of the standard tables, [**26**, **34**] and which is a special case of a more general result [**26**, Eq. 7.4.4.11 or 7.4.4.3] follows:

$$ {}_3F_2\left({1,n+1,\beta \atop n+2,\alpha}\Big|1\right)=\frac{(n+1)\Gamma(\alpha)}{\Gamma(\beta)}\left(\frac{\Gamma(\beta-n-1)(\psi(\alpha-n-1)-\psi(\alpha-\beta))}{\Gamma(\alpha-n-1)} - \sum_{l=0}^{n-1}\frac{\Gamma(\beta+l-n)}{\Gamma(\alpha+l-n)(l+1)}\right) . \tag{12}$$

**Proof:** Write the ${}_3F_2$ as a convergent sum with $\Re(\alpha)>\Re(\beta)$, change the lower summation limit from "0" to "-n", subtracting equal terms, to get:

$$\begin{aligned} {}_3F_2\left({1,n+1,\beta \atop n+2,\alpha}\Big|1\right) &= \frac{(n+1)\Gamma(\alpha)}{\Gamma(\beta)}\sum_{l=0}\frac{\Gamma(\beta+l)}{\Gamma(\alpha+l)(l+n+1)} \\ &= \frac{(n+1)\Gamma(\alpha)}{\Gamma(\beta)}\left(\sum_{l=0}\frac{\Gamma(\beta+l-n)}{\Gamma(\alpha+l-n)(l+1)}-\sum_{l=0}^{n-1}\frac{\Gamma(\beta+l-n)}{\Gamma(\alpha+l-n)(l+1)}\right) . \end{aligned} \tag{13}$$

The first infinite sum in the second equality is recognized as a ${}_3F_2$ with known summation [**18**, Eq. 3.13(43)]. **(12)** follows immediately, for all values of α and β by the principle of analytic continuation. §

**Lemma 2.2**

Another useful result needed later (see **(63)**) that does not appear in the standard tables, (however, see [**22**]) but which is a special case of a known result [**26**, Eq. 7.4.4.11] is

$$\begin{aligned} {}_3F_2\left({a,b,c \atop n+b,c+1,}\Big|1\right) &= \frac{(b)_n\Gamma(c+1)\Gamma(1-a)}{(b-c)_n\Gamma(c+1-a)}+ \\ & c\Gamma(b+n)\Gamma(c-b+1-n)\sum_{\ell=0}^{n-1}\frac{\Gamma(n-\ell-a)(-1)^{\ell}}{\Gamma(b+n-a-\ell)\Gamma(n-\ell)\Gamma(c-b-n+2+\ell)}. \end{aligned} \tag{14}$$

**Proof:** With $\Re(n+1)>\Re(\alpha)$, apply the well-known result [**18**, Eq. 3.13(38)] n-times repeatedly and use Gauss' formula for a ${}_2F_1$ of unit argument. §

**Corollary:** Evaluate the limit a=1 in **(9**:

$$
{}_3F_2\left({1,b,c \atop n+b,c+1,}\Big|1\right)=\frac{c(\psi(b)-\psi(c))\Gamma(n+b)\Gamma(b-c)}{\Gamma(b)\Gamma(b-c+n)}
+c\Gamma(b+n)\Gamma(c-b+1-n)\sum_{\ell=0}^{n-2}\frac{\Gamma(n-\ell-1)(-1)^{\ell}}{\Gamma(b+n-1-\ell)\Gamma(n-\ell)\Gamma(c-b-n+2+\ell)}, \qquad \textbf{(15)}
$$

reducing to [**26**, Eq. 7.4.4(33)] when n=1.

**Theorem 2.1**

The following terminating sum only appears symbolically in the literature as a limiting case of a very complicated result [**26**, Eq. 7.10.2(9)], or for special values of "*c*" (e.g. [**29**, Eq. (3.4)]):

$$
{}_4F_3\left({1,1,1,-k \atop 2,2,1+c}\Big|1\right)=\frac{c}{(k+1)}\left(\sum_{l=0}^{k}\frac{\psi(c+l+1)}{(1+l)}-\psi(c)\left(\psi(k+2)-\psi(1)\right)\right) \qquad \textbf{(16a)}
$$

$$
=\frac{c}{(k+1)}\left(\psi(k+2)\left(\psi(c+k+1)-\psi(c)\right)-\sum_{l=0}^{k}\frac{\psi(l+1)}{(c+l)}\right). \qquad \textbf{(16b)}
$$

The two forms of **(16)** are obtained from one another using **(21)** below. The derivation depends on results to be obtained in Section 4, but this result is listed here for consistency of presentation. The proof is outlined in Appendix A.

**Theorem 2.2**

Similarly, the following, which only appears in the tables as an exceptional, limiting case of a very complicated result [**26**, Eq. 7.10.2(2)], is easily obtained by rewriting **(C.20)**.

$$
\frac{1}{q}\,{}_4F_3\left({1,1,1,1 \atop 2,2,1+q}\Big|1\right)=\tfrac{1}{2}\psi^{(2)}(q)-[\psi'(q)+\frac{\pi^2}{6}][(\gamma+\psi(q)]+2\zeta(3)
+\Gamma(q)\sum_{l=0}\frac{(-1)^l\psi'(l+2)}{\Gamma(q-l-1)\Gamma(l+2)(l+1)}\quad ,\Re(1+q)>0\,. \qquad \textbf{(17)}
$$

If $q\to n,\ n\geq 1$, the infinite series terminates at $l=n-2$.

**Theorem 2.3**

Likewise, the following results do not appear in the literature.

| | |
|---|---|
| $\dfrac{1}{(\frac{1}{2}+n)}{}_4F_3\left({1,1,1,1 \atop 2,2,3/2+n}\,\middle\vert\,1\right)=\frac{7}{2}\zeta(3)+\frac{1}{2}\psi^{(2)}(\tfrac{1}{2}+n)-\dfrac{\pi^2}{2}\psi(\tfrac{1}{2}+n)-\gamma\psi'(\tfrac{1}{2}+n)$ $+\sum_{l=0}^{n-1}\left(\dfrac{\psi(1+l)}{(l-n+\frac{1}{2})^2}-\psi'(1+l)\left(\dfrac{1}{(l+\frac{1}{2})}+\dfrac{1}{(l-n+\frac{1}{2})}\right)\right)$ | **(18)** |
| $\dfrac{1}{(\frac{1}{2}-n)}{}_4F_3\left({1,1,1,1 \atop 2,2,3/2-n}\,\middle\vert\,1\right)=\frac{7}{2}\zeta(3)+\frac{1}{2}\psi^{(2)}(\tfrac{1}{2}+n)-\pi^2\left[\psi(\tfrac{1}{2}+n)+\tfrac{3}{2}\gamma+\log 2\right]$ $+\psi'(\tfrac{1}{2}+n)\left[2\gamma+2\log 2+\psi(n+\tfrac{1}{2})\right]$ $-\sum_{l=0}^{n-1}\left(\dfrac{\psi(1+l)}{(l+\frac{1}{2})^2}-\psi'(1+l)\left(\dfrac{1}{(l+\frac{1}{2})}+\dfrac{1}{(l-n+\frac{1}{2})}\right)\right)$ | **(19)** |
| ${}_4F_3\left({1,1,1,1 \atop 2,2,\,1/2}\,\middle\vert\,1\right)=\frac{7}{4}\zeta(3)+\dfrac{\pi^2}{2}(1-\log 2)$ | **(20a)** |
| ${}_4F_3\left({1,1,1,1 \atop 2,2,\,3/2}\,\middle\vert\,1\right)=-\frac{7}{4}\zeta(3)+\dfrac{\pi^2}{2}\log 2$ | **(20b)** |

**Proof:**

Set $q=\frac{1}{2}\pm n$ in (**C.20**) employing **(3)** and (**C.23**). §

In **(19)**, since the series representation of the left-hand side does not converge for $n \geq 2$, the right-hand side represents the left-hand side* in the sense of analytic continuation, because the variable "*q*" is continuous in (**C.20**) (whose left-hand side is a valid representation of the hypergeometric function for this range of the variable "*n*"). **(20a)** is a simple restatement and simplification of **(19)** using $n$=1. The case n=0, **(20b)**, corresponds to the special case $a=b\to 0$ of [**15**, Eq.7].

## 3. Finite sums involving Digamma functions

**Theorem 3.1**

An important duality exists between two finite sums of digamma functions (see also (**B.6**)):

$$\sum_{l=0}^{k}\frac{\psi(b+l)}{c+l}=\psi(c+k+1)\psi(b+k+1)-\psi(b)\psi(c)-\sum_{l=0}^{k}\frac{\psi(c+l+1)}{b+l}\,. \qquad \textbf{(21)}$$

**Proof:** Expand, rearrange, and with reference to **(3)** and **(5)**, find:

* i.e. ${}_4F_3\left({1,1,1,1 \atop 2,2,3/2-n}\,\middle\vert\,1\right)$ means $\lim_{q\to n}\,{}_4F_3\left({1,1,1,1 \atop 2,2,3/2-q}\,\middle\vert\,1\right),\quad n\geq 0$ .

$$\begin{aligned}\sum_{l=0}^{k}\frac{\psi(b+l)}{c+l} &= \frac{\psi(b)}{c}+\frac{\psi(b+1)}{c+1}+\frac{\psi(b+2)}{c+2}+\cdots+\frac{\psi(b+k)}{c+k}\\ &=(\psi(b)+\tfrac{1}{b}+\tfrac{1}{b+1}+\cdots\tfrac{1}{b+k-2}+\tfrac{1}{b+k-1})/(c+k)+\\ &\quad(\psi(b)+\tfrac{1}{b}+\tfrac{1}{b+1}+\cdots\tfrac{1}{b+k-2}\qquad)/(c+k-1)+\cdots+\cdots+\\ &\quad(\psi(b)\qquad\qquad\qquad\qquad)/(c)\\ &=\psi(b)(\psi(c+k+1)-\psi(c))+\sum_{l=0}^{k-1}\frac{\psi(c+k+1)-\psi(c+1+l)}{b+l}\end{aligned} \tag{22}$$

from which **(21)** follows immediately. §

Using **(3)**, the sums on either side of the equality can be reduced to one another when *"b"* and *"c"* are separated by integers, showing that a closed sum exists for those cases (cf. **(26)** and **(B.3)**). Another, equivalent form of **(21)** can be easily found by rearrangement:

$$\begin{aligned}\sum_{l=0}^{k}\left(\frac{\psi(b+l)}{c+l}+\frac{\psi(c+l)}{b+l}\right) &= \psi(c+k+1)\psi(b+k+1)-\psi(b)\psi(c)\\ &\quad+\frac{\psi(b+k+1)-\psi(c+k+1)}{b-c}-\frac{\psi(b)-\psi(c)}{b-c}.\end{aligned} \tag{23}$$

A useful [**19**] special case is given below:

$$\sum_{l=0}^{k}\frac{\psi(l+\frac{1}{2})}{l+1} = \psi(k+\tfrac{1}{2})\psi(k+2)-\psi(1)\psi(\tfrac{1}{2})-\sum_{l=1}^{k}\frac{\psi(l+1)}{l-\frac{1}{2}} \tag{24}$$

<u>**Corollary**</u>

Apply **(3)** to **(21)** and simplify to obtain an equivalent form in terms of generalized harmonic numbers [**6**, **35**]:

$$\sum_{l=0}^{k}\frac{1}{l+1}H_l(\tfrac{1}{2})+\sum_{l=0}^{k}\frac{1}{l-\frac{1}{2}}H_l(1) = H_k(\tfrac{1}{2})H_{k+1}(1) \tag{25}$$

Also see [**30**, Theorem 7, Example 2].

<u>**Theorem 3.2**</u> (See also **(B.9)**)

$$\sum_{l=0}^{k}\frac{\psi(b+l)}{(b+l)} = \tfrac{1}{2}[\psi'(b+1+k)-\psi'(b)+\psi(b+1+k)^2-\psi(b)^2] \tag{26}$$

<u>**Proof:**</u> Start from [**12**, Eq.7.1.1]

$$\sum_{l=0}^{k}\frac{\Gamma(b+l))}{\Gamma(a+l)}=(\frac{\Gamma(b+1+k)}{\Gamma(a+k)}-\frac{\Gamma(b)}{\Gamma(a-1)})/(b-a+1) \qquad \textbf{(27)}$$

and operate on both sides with $\frac{\partial}{\partial b}$ giving an interim result worth noting:

| | |
|---|---|
| $$\sum_{l=0}^{k}\frac{\Gamma(b+l)\psi(b+l)}{\Gamma(a+l)}=(\frac{\Gamma(b+1+k)\psi(b+1+k)}{\Gamma(a+k)}-\frac{\Gamma(b)\psi(b)}{\Gamma(a-1)})/(b-a+1) -(\frac{\Gamma(b+1+k)}{\Gamma(a+k)}-\frac{\Gamma(b)}{\Gamma(a-1)})/(b-a+1)^2 \;.$$ | **(28)** |

Then evaluate the limit $a \to b+1$ to get **(26)**. Alternatively, take the limit $c=b$ in **(23)**. §

In applications, the cases $b=1$ and $b=1/2$, both of which are easily obtained from **(26)**, arise.* Notably, the case $b=1$ extends an unfinished sum† and integral from the standard tables [**12**, Eqs. (55.2.1) and (5.13.17)] as follows:

$$-\int_0^1\frac{1-t^k}{1-t}\log(1-t)dt=\sum_{l=1}^{k}\frac{1}{l}\sum_{m=1}^{l}\frac{1}{m}=\tfrac{1}{2}[(\psi(1+k)-\psi(1))^2+\psi'(1)-\psi'(1+k)]$$
$$=\frac{1}{2}\left[\sum_{l=0}^{k-1}\frac{1}{(l+1)^2}+\left(\sum_{l=0}^{k-1}\frac{1}{(l+1)}\right)^2\right] \qquad \textbf{(29)}$$

specializing‡ a recently quoted general result for non-integral values of "*k*" [**3**, Eqs.(3.6) and (3.9)]:

$$\int_0^1\frac{1-t^x}{1-t}\log(1-t)dt=\sum_{l=2}^{\infty}(-x)^{l-1}\sum_{m=1}^{\infty}\frac{\psi(m)+\gamma}{m^l}$$
$$=-\tfrac{1}{2}[(\psi(1+x)-\psi(1))^2+\psi'(1)-\psi'(1+x)] \qquad \textbf{(30)}$$

## 4. Infinite Series

For completeness, a well-known sum, recently revisited [**27**, Eq. (2.7)], is quoted. If $\Re(c-b)>0$,

| | |
|---|---|
| $$\sum_{l=1}^{\infty}\frac{\Gamma(b+l)}{\Gamma(c+l)}(\frac{\psi(b+l)-\psi(b+1)}{l})=\frac{\Gamma(b)}{\Gamma(c)}\{\psi'(c-b)-\frac{[\psi(c)-\psi(c-b)]}{b}\}$$ | **(31)** |

* Curiously, the computer code Maple 6 (and Maple 2016) found the b=1 version of (the new result) **(26)**, but missed all other possibilities tested. Mathematica [**34**,Version 4.1 (~2004)] missed all possibilities tested, but Mathematica [**10**, version(2016)] correctly found this result.

† reproducing a known result [**1**, Lemma 1] in this special case.

‡ Note the relationship between the inner sum in **(30)** and **(67)**.

Set c=b+1, and change the notation slightly, giving

$$\sum_{l=0}^{\infty}\frac{\psi(b+l+1)}{(l+1)(b+l+1)}=\frac{1}{b}\{\frac{\pi^2}{6}+\psi(b)[\psi(b+1)-\psi(1)]\}\,. \tag{32}$$

In the limit *b→0* this reduces to a very well-known result [**9**,**24**, **34**]

$$\sum_{l=0}^{\infty}\frac{\psi(l+1)}{(l+1)^2}=\zeta(3)-\frac{\pi^2\gamma}{6} \tag{33}$$

**Theorem 4.1**

Consider a useful variant of **(32)**:

| | |
|---|---|
| $$\sum_{l=0}\frac{\psi(l+1)}{(l+p)(l+q)}=\frac{1}{2(q-p)}\{\psi'(p)-\psi(p)^2-\psi'(q)+\psi(q)^2\}$$ | **(34)** |

**Proof:** Temporarily, assume $\Re(\beta)<1$, and consider

$$\begin{aligned}\sum_{l=0}\frac{\Gamma(l+\beta)}{\Gamma(l+1)}\frac{1}{(l+p)(l+q)}&=\frac{1}{(q-p)}\sum_{l=0}\frac{\Gamma(l+\beta)}{\Gamma(l+1)}[\frac{1}{(l+p)}-\frac{1}{(l+q)}]\\&=\frac{1}{(q-p)}[\frac{\Gamma(\beta)\Gamma(p)}{\Gamma(p+1)}{}_2F_1(\begin{matrix}\beta,p\\p+1\end{matrix}1)-\frac{\Gamma(\beta)\Gamma(q)}{\Gamma(q+1)}{}_2F_1(\begin{matrix}\beta,q\\q+1\end{matrix}1)]\\&=\frac{\Gamma(\beta)\Gamma(1-\beta)}{(q-p)}(\frac{\Gamma(p)}{\Gamma(p+1-\beta)}-\frac{\Gamma(q)}{\Gamma(q+1-\beta)})\end{aligned} \tag{35}$$

using Gauss' summation for a $_2F_1$ of unit argument. The condition on $\beta$ may now be relaxed. Operate on both sides of **(35)** with $\lim_{\beta\to 1}\frac{\partial}{\partial\beta}$, respectively reducing the left- and right-hand sides to the corresponding sides of **(34)**. §

The sums (enclosed in square brackets) in these results diverge if they are considered individually when $\Re(\beta)\geq 1$, so that a careful adherence to the ordering of limits, and the use of the principle of analytic continuation is needed.

**Corollary**

The following generalizes **(33)**

| | |
|---|---|
| $$\sum_{l=0}\frac{\psi(l+1)}{(l+q)^2}=\psi'(q)\psi(q)-\tfrac{1}{2}\psi^{(2)}(q)$$ | **(36)** |

**Proof:** Evaluate the limit $p\to q$ in **(34)**. § See also **(66)** and ( .

**Corollary**

| | |
|---|---|
| $$\sum_{l=0}\frac{\psi(l+m+1)}{(l+p)(l+q)}=\frac{1}{(q-p)}\{(\psi'(p)-\psi(p)^2-\psi'(q)+\psi(q)^2)/2$$ $$+\psi(p)^2-\psi(q)^2-\psi(p)\psi(p-m)+\psi(q)\psi(q-m)\}-\sum_{l=0}^{m-1}\frac{\psi(l+1)}{(p-l-1)(q-l-1)}$$ | **(37)** |

**Proof:** Using **(3)** and **(6)** assume $p,q\neq 1,2,\cdots m$ and split the following sum using partial fractions

$$\sum_{k=0}^{m-1}\sum_{l=0}\frac{1}{(l+p)(l+q)(l+m-k)}$$
$$=\frac{\psi(p)}{(q-p)}\sum_{l=0}^{m-1}\frac{1}{p-m+l}-\frac{\psi(q)}{(q-p)}\sum_{l=0}^{m-1}\frac{1}{q-m+l}-\sum_{l=0}^{m-1}\frac{\psi(m-l)}{(q-m+l)(p-m+l)} \quad \textbf{(38)}$$

The first two sums are easily evaluated; reverse the third right-side sum, apply **(34)** to the finite sum on the left and eventually arrive at **(37)**. §

The conditions on $p$ and $q$ can be relaxed by taking the appropriate limits. A limit must also be evaluated when $p$ and $q$ are separated by integers. In particular, many results, laboriously and individually derived elsewhere [**24**] are special cases of **(37)**. Also, if $p$ and $q$ are integers, the finite sum in **(37)** can be evaluated in closed form using **(26)** – also see **(B.3)**, (**B.16**), (**C.6)**) and **(C.7)** as well as Xu [**36**]).

**Theorem 4.2**

| | |
|---|---|
| $$\sum_{l=0}\frac{[\psi(\alpha+q+l)-\psi(\alpha+l)]}{l+1}=\rho(\alpha,q) \qquad \text{where}$$ $$\rho(\alpha,q)\equiv\frac{1}{2}(\psi(\alpha+q-1)^2-\psi(\alpha-1)^2-\psi'(\alpha+q-1)+\psi'(\alpha-1))$$ $$-\psi(1)(\psi(\alpha+q-1)-\psi(\alpha-1))$$ | **(39)** |

**Proof:** Consider the following sum with $\Re(\alpha)>\Re(\beta)$:

$$\sum_{l=0}[\frac{\Gamma(\beta+q+l)}{\Gamma(\alpha+q+l)}-\frac{\Gamma(\beta+l)}{\Gamma(\alpha+l)}]\frac{1}{(l+1)}=\frac{\Gamma(\beta+q)}{\Gamma(\alpha+q)}\,{}_3F_2(\begin{matrix}\beta+q,1,1\\ \alpha+q,2\end{matrix}|1)-\frac{\Gamma(b)}{\Gamma(a)}\,{}_3F_2(\begin{matrix}\beta,1,1\\ \alpha,2\end{matrix}|1)$$
$$=\frac{\Gamma(\beta+x-1)}{\Gamma(\alpha+x-1)}\psi(\alpha+q-1)-\frac{\Gamma(\beta-1)}{\Gamma(\alpha-1)}\psi(\alpha-1)+\psi(\alpha-b)\Big(\frac{\Gamma(\beta-1)}{\Gamma(\alpha-1)}-\frac{\Gamma(\beta+q-1)}{\Gamma(\alpha+q-1)}\Big) \quad \textbf{(40)}$$

using a well-known result [**18**, Eq. 3.13(43)] to sum the hypergeometric functions. Relax the constraint and operate on both sides of **(40** with $\lim_{\beta\to\alpha}\frac{\partial}{\partial\beta}$ so that the right- and left-hand sides become the corresponding sides of **(39)**. §

As before, the two terms in **(40)** correspond to divergent series if written individually; this is indicated by the use of square brackets as a reminder, and again new results can be obtained by operating with $\frac{\partial}{\partial\alpha}$ and/or $\frac{\partial}{\partial q}$ and taking appropriate limits. For example, the limit $\alpha=1$ in **(39)** yields

| | |
|---|---|
| $\sum_{l=0}\frac{\psi(q+l+1)-\psi(l+1)}{(l+1)}=\left((\psi(q)-\psi(1))^2+\psi'(\tfrac{1}{2})-\psi'(q)\right)/2$ | **(41)** |

and operating with $\frac{\partial}{\partial q}$ on **(41)** with some reordering gives

| | |
|---|---|
| $\sum_{l=0}\frac{\psi'(q+l)}{(l+1)}=-\frac{1}{2}\psi^{(2)}(q-1)+\psi(q-1)\psi'(q-1)-\psi(1)\psi'(q-1)\,.$ | **(42)** |

Take the limit q=1 to recover the special case

$$\sum_{l=0}\frac{\psi'(1+l)}{(l+1)}=2\zeta(3)\,, \tag{43}$$

corresponding to a known result [**24**], [**12**, Eq. 55.9.7).

**Theorem 4.3**
Eq. **(39)** can be generalized further.

| | |
|---|---|
| $\sum_{l=0}\frac{[\psi(\alpha+q+l)-\psi(\alpha+l)]}{l+m+1}=\rho(\alpha-m,q)-\sum_{l=0}^{m-1}\frac{\psi(\alpha+q+l-m)-\psi(\alpha+l-m)}{l+1}$ | **(44)** |

**Proof:** Add and subtract terms corresponding to $l=-m,\cdots,-1$ in **(39)** and re-order the resulting series. §

**Corollary**: **(44)** reduces to the following useful result by setting $\alpha=m+\frac{1}{2}$ and $q=\frac{1}{2}$.

$$\sum_{l=0}\frac{[\psi(l+m+1)-\psi(l+m+\frac{1}{2})]}{(l+m+1)}=\tfrac{1}{2}\left(\psi'(1)-\psi(\tfrac{1}{2})^2-\psi'(m+1)-\psi(m+1)^2\right)+2\left(\psi(1)-\psi(\tfrac{1}{2})\right)+\psi(\tfrac{1}{2})\psi(1)+\sum_{l=0}^{m-1}\frac{\psi(l+\frac{1}{2})}{(l+1)} \qquad \textbf{(45)}$$

See also (**24**).

**Theorem 4.4**

For completeness' sake, a known variation of **(31)** [**27**, Eq. (2.14)] is quoted, without proof.

$$\sum_{l=0}\frac{\Gamma(q+l)}{\Gamma(p+l)}\psi(l+1)=\frac{\Gamma(q)\Gamma(p-q-1)}{\Gamma(p-1)\Gamma(p-q)}\left(\psi(1)-\psi(p-q-1)+\psi(p-1)\right) \qquad \textbf{(46)}$$

Using methods similar to those already employed (e.g. Theorem 4.2), other variants of **(31)** and **(46)** can be obtained:

$$\sum_{l=0}\frac{\Gamma(q+l)\psi(p+l)}{\Gamma(p+l)(l+1)}=\frac{\Gamma(q-1)}{\Gamma(p-1)}[\psi'(p-q)-\psi'(p-1)+\psi(p-1)[\psi(p-1)-\psi(p-q)]] \qquad \textbf{(47)}$$

and

$$\sum_{l=0}\frac{\Gamma(q+l)\psi(1+l)}{\Gamma(p+l)(1+l)}=\frac{\Gamma(q-1)}{2\Gamma(p-1)}[[\psi(p-1)-\psi(p-q)]^2-2\gamma[\psi(p-1)-\psi(p-q)]+\psi'(p-q)-\psi'(p-1)] \qquad \textbf{(48)}$$

**Theorem 4.5**

$$\sum_{l=0}\frac{\Gamma(c+l)}{\Gamma(p+l)}\psi(f+l)=\frac{\Gamma(c)}{(p-c-1)}\left\{\frac{(\psi(f-1)+\frac{1}{(p-c-1)})}{\Gamma(p-1)}-\frac{\Gamma(f)\Gamma(1-f)\Gamma(1-c)}{\Gamma(f-c)\Gamma(p-f)}-\frac{\Gamma(2-f)\Gamma(p-c)}{\Gamma(p-f)\Gamma(f-c)}\sum_{l=0}\frac{\Gamma(f-c+l)(-1)^l}{\Gamma(2-c+l)\Gamma(l+1)\Gamma(p-1-l)(1-c+l)}\right\} \qquad \textbf{(49)}$$

**Proof:** Consider the sum

$$\sum_{l=0}\frac{\Gamma(c+l)\Gamma(\beta+l)}{\Gamma(p+l)\Gamma(f+l)}=\frac{\Gamma(c)\Gamma(\beta)}{\Gamma(p)\Gamma(f)}\,{}_3F_2\left(\begin{matrix}1,\beta,c\\p,f\end{matrix}\,\middle|\,1\right) \qquad \textbf{(50)}$$

Convert the hypergeometric ${}_3F_2$ into another ${}_3F_2$ having the property that it terminates if $p=m$, ($m$>1) according to [**18**, section 3.13.3]

$$ {}_3F_2\left({1,\beta,c \atop p,f}\,\middle|\,1\right) = \frac{\Gamma(p+f-1-\beta-c)\Gamma(p)\Gamma(1-\beta)}{\Gamma(p-\beta)}\left\{\frac{\Gamma(f)\Gamma(1-c)}{\Gamma(f-c)\Gamma(p-c)\Gamma(f-\beta)} - \frac{(f-1)}{(1-c)\Gamma(f-\beta-c+1)\Gamma(p-1)}\,{}_3F_2\left({f-c,2-p,1-c \atop f-c-\beta+1,2-c}\,\middle|\,1\right)\right\} \qquad \textbf{(51)} $$

Operate on both sides with $\lim\limits_{\beta \to f} \frac{\partial}{\partial \beta}$ and simplify. §

This derivation of **(49)** is valid for all continuous values of $p$, although it is not particularly useful unless the infinite series (equivalent to a ${}_3F_2$) can be summed analytically (e.g. $p = m$, an integer, in which case the infinite series terminates). If $f$=1, it can be shown, by summing the infinite series [**26**, Eq.7.4.4.(49)] and taking the appropriate limit, that **(49)** reduces to **(46)** as a special case. To obtain other useful results from **(49)**, operate with $\lim\limits_{f \to 1} \frac{\partial}{\partial f}$. Then set $p = 2$ giving

| | |
|---|---|
| $\displaystyle\sum_{l=0} \frac{\Gamma(l+q)\psi'(l+1)}{\Gamma(l+2)} = \Gamma(q-1)\left(\tfrac{1}{2}\psi'(1-q) - \tfrac{1}{2}\psi(1-q)^2 - \tfrac{1}{4}\pi^2 - \tfrac{1}{2}\gamma^2 - \gamma\psi(1-q)\right),$ | **(52)** |

reducing to **(43)** when $q$=1.

Alternatively, operate on **(49)** with $\lim\limits_{p \to c+1} \lim\limits_{f \to 1} \frac{\partial}{\partial f}$, and from the requirement that the coefficient of $(p-c-1)^{-1}$ on the right-hand side must vanish (series converges), find

| | |
|---|---|
| $\displaystyle\sum_{l=0} \frac{\Gamma(l+q)\psi(l+q)}{\Gamma(l+1)(q+l)^2} = \frac{\pi\Gamma(q)}{\sin \pi q}\left[(\psi(q)+\gamma)(\gamma + \pi\cot\pi q) - \frac{\pi^2}{6}\right],$ | **(53)** |

reducing to **(33)** when $q \to 1$. See also **(C.19)** to (**C.23**).

**Lemma 4.6**

The following sums the difference of two divergent series (see also **(C.43)** and **(C.44))**.

| | |
|---|---|
| $\displaystyle\sum_{l=0} \frac{\psi(l+\frac{1}{2}) - \psi(l+1)}{l+\frac{1}{2}} = -\frac{\pi^2}{3}$ | **(54)** |

**Proof:** Consider the following sum (with $\Re(\beta) < 0$):

$$\sum_{l=0}\frac{\Gamma(l+\beta-\frac{1}{2})\Gamma(l+1)-\Gamma(l+\beta)\Gamma(l+\frac{1}{2})}{\Gamma(l+\frac{3}{2})\Gamma(l+1)}$$
$$=\frac{\Gamma(\beta-\frac{1}{2})}{\Gamma(\frac{3}{2})}{}_2F_1\left({1,\beta-1/2 \atop 3/2}\middle|1\right)-\frac{\Gamma(\beta)\Gamma(\frac{1}{2})}{\Gamma(\frac{3}{2})}{}_2F_1\left({\beta,1/2 \atop 3/2}\middle|1\right). \tag{55}$$

As before, sum the hypergeometric series, relax the limitation on $\beta$ and, with deference to the principle of analytic continuation, operate on both sides with $\lim_{\beta\to 1}\frac{\partial}{\partial\beta}$. §

**Theorem 4.7**

The following generalizes **(54)** (also see **(C.9)** etc.):

| | |
|---|---|
| $$\sum_{l=0}\left[\frac{\psi(l+k+\frac{1}{2})-\psi(l+k+1)}{l+\frac{1}{2}}\right]=\psi(k)\psi(k+\tfrac{1}{2})-(\psi(\tfrac{1}{2})+4)(\psi(k)-\psi(k+\tfrac{1}{2})) +\frac{\pi^2}{6}-\psi(\tfrac{1}{2})^2-8-4\psi(\tfrac{1}{2})-2\sum_{l=1}^{k-1}\frac{\psi(l)}{(l+\frac{1}{2})}$$ | **(56)** |

**Proof:** Let

$$V_k=\sum_{l=0}\frac{\psi(l+k+\frac{1}{2})-\psi(l+k+1)}{l+\frac{1}{2}} \tag{57}$$

Then

$$V_k=V_{k-1}+\sum_{l=0}\left(\frac{1}{l+k-\frac{1}{2}}-\frac{1}{l+k}\right)\left(\frac{1}{l+\frac{1}{2}}\right)=V_{k-1}+\frac{\psi(\frac{1}{2})-\psi(k)}{k-\frac{1}{2}}-\frac{\psi(\frac{1}{2})-\psi(k-\frac{1}{2})}{k-1} \tag{58}$$

Apply **(58)** repeatedly (k times), giving

$$V_k=V_0+\sum_{l=1}^{k}\left(\frac{\psi(\frac{1}{2})-\psi(l)}{l-\frac{1}{2}}-\frac{\psi(\frac{1}{2})-\psi(l-\frac{1}{2})}{l-1}\right)\qquad with\quad V_0=-\frac{\pi^2}{3} \tag{59}$$

The second term enclosed by brackets in **(59)** equals $\pi^2/2\ if\ l=1$. Use **(54)** and its analogues, and, after some simplification, including the use of **(21)**, find **(56)**. § See also **(45)**.

**Theorem 4.8**

| | |
|---|---|
| $$\sum_{l=0}\left[\frac{\psi(a+l)}{(a+l)}-\frac{\psi(a+l+q)}{(a+l+q)}\right]=\left(\{\psi(a+q)^2+\psi'(a+q)\}-\{\psi(a)^2+\psi'(a)\}\right)/2$$ | **(60)** |

**Proof:** Since this is effectively the (finite) difference between diverging sums, consider the following two converging sums, subject to the condition $\Re(a)>\Re(\beta)$:

$$[\sum_{l=0} \frac{\Gamma(\beta+l)}{\Gamma(a+l+1))} - \sum_{l=0} \frac{\Gamma(\beta+l+q)}{\Gamma(a+l+q+1)}]$$
$$= \frac{\Gamma(\beta)}{\Gamma(a+1)} {}_2F_1\left(\begin{matrix}\beta,1\\ a+1\end{matrix}\,|\,1\right) - \frac{\Gamma(\beta+q)}{\Gamma(a+q+1)} {}_2F_1\left(\begin{matrix}\beta+q,1\\ a+q+1\end{matrix}\,|\,1\right) \qquad \textbf{(61)}$$

As in previous cases, evaluate the two hypergeometric functions independently, simplify, and operate with $\frac{\partial}{\partial\beta}$. Then, unite the two sums into one (convergent) sum, and evaluate the limit $\beta \to a$. §

Notice that, to within a constant, each of the terms of **(60)** that are enclosed in braces ({ }) regularize corresponding divergent series on the left. Also, see **(C.8)** and **(C.14)**.

**Corollary**

| | |
|---|---|
| $\sum_{l=0} \frac{\psi'(q+l)}{(q+l)} - \frac{\psi(q+l)}{(q+l)^2} = -\frac{1}{2}\psi^{(2)}(q) - \psi(q)\psi'(q)$ | **(62)** |

**Proof**: Operate with $\lim_{a\to 0} \frac{\partial}{\partial q}$ on **(60)**; alternatively, operate with $\lim_{b\to q} \lim_{k\to\infty} \frac{\partial}{\partial b}$ on **(26)**. §

**Theorem 4.9**

| | |
|---|---|
| $\sum_{l=0}[\frac{\psi(c+l)}{(n+l+1)} - \frac{\psi(n+l+1)}{(c+l)}] = \frac{\pi^2}{3} - \psi(n+1)(\psi(c) + \frac{1}{(c-n-1)})$ $-\psi'(c-n-1) + \psi(c-n-1)^2 + \frac{\psi(c)}{(c-n-1)} + 2\sum_{l=0}^{n} \frac{\psi(l+1)}{(c+l-n-1)}$ | **(63)** |

**Proof:** Consider the following sum with $\Re(c-\beta) > 0$. With this restriction, the sum may be split into two and identified:

$$\sum_{l=0}[\frac{\Gamma(\beta+l)}{\Gamma(c+l)(n+1+l)} - \frac{\Gamma(\beta+l+n+1-c)}{\Gamma(n+1+l)(c+l)}]$$
$$= \frac{\Gamma(\beta)}{(n+1)\Gamma(c)} {}_3F_2\left(\begin{matrix}1,n+1,\beta\\ c,n+2\end{matrix}\,|\,1\right) - \frac{\Gamma(\beta+n+1-c)}{(c)\Gamma(n+1)} {}_3F_2\left(\begin{matrix}1,c,\beta+n+1-c\\ n+1,c+1\end{matrix}\,|\,1\right) \qquad \textbf{(64)}$$

Using **(12)** and **(9**, each of the hypergeometric functions can be summed, and the restriction relaxed by the principle of analytic continuation. Since the combined sums converge, operate on both sides with $\lim_{\beta\to c} \frac{\partial}{\partial\beta}$ and sum a terminating $_3F_2$, giving

$$\sum_{l=0}[\frac{\psi(c+l)}{(n+l+1)}-\frac{\psi(n+l+1)}{(c+l)}]=\frac{\pi^2}{3}+\gamma\psi(c-1-n)+\psi(n+1)(\psi(c)-\psi(c-n))$$
$$-\psi(1,c-n)+\psi(c-n)^2-\frac{\psi(c-n)}{(c-n-1)}-\sum_{l=0}^{n-1}\frac{\psi(c+l-n)}{(l+1)}-\frac{n}{(c-n)}{}_4F_3\left(\begin{matrix}1-n,1,1,1\\c+1-n,2,2\end{matrix}\Big|1\right) \quad \textbf{(65)}$$

From either version of **(16)**, the $_4F_3$ can be identified, and Eq. **(63)** emerges after considerable simplification. § See also (**B.16**).

In the limit $c \to n+1$, both sides of **(63)** and **(65)** vanish. This has been verified and a known, closed-form for a special case of the resulting $_4F_3$ has been retrieved.

**Theorem 4.10** (generalized Euler sums [**10**] and polylogarithms [**16**]. See also **(C.14)**.

| | |
|---|---|
| $\sum_{l=0}\frac{\psi(l+1)}{(l+q)^n}=\frac{(-1)^n}{(n-1)!}\left[\sum_{k=0}^{n-2}\psi^{(1+k)}(q)\psi^{(n-k-2)}(q)\binom{n-2}{k}-\frac{1}{2}\psi^{(n)}(q)\right]\quad ,n>2$ | **(66)** |

**Proof**: Operate on **(36)** with $\frac{\partial^{n-2}}{dq^{n-2}}$. §

Using $\psi^{(n)}(1)=(-1)^{n+1}n!\zeta(n+1)$ and setting q=1 gives a simple derivation for an equivalent form of Euler's (linear) sum:

$$\sum_{l=0}\frac{\psi(l+2)-\psi(1)}{(l+1)^n}=(1+\frac{n}{2})\zeta(n+1)-\frac{1}{(n-1)}\sum_{k=1}^{n-2}k\zeta(k+1)\zeta(n-k) \quad \textbf{(67)}$$

which, when compared to the form usually quoted in the literature [**1**, Eq. 20], [**10**, Theorem 2.2], yields the interesting relation

$$\sum_{k=1}^{n-2}\zeta(k+1)\zeta(n-k)=\frac{2}{(n-1)}\sum_{k=1}^{n-2}k\,\zeta(k+1)\zeta(n-k), \quad \textbf{(67a)}$$

a result that is true for any indexed function (proof by reversing the right-hand sum). Xu **[36, Eqs. (1.27)** and **(1.28)** – incorrect – add a term $1/a^{s+1}$ to the latter) has recently obtained an equivalent form of **(66)** as have Sofo and Cvijovoć **[28, Theorem 1].** Special cases of **(66)** with combinations of $q = 0, \frac{1}{2}, 1$ and $n = 2, 3$ can be found in Appendix E; for example, set n=3 in **(66)** to retrieve the well-known harmonic result [**9**, Eq. (17)]:

$$\sum_{l=0}\frac{\psi(l+1)-\psi(1)}{(l+1)^3}=\frac{\pi^4}{360}. \quad \textbf{(68)}$$

Finally, by taking the limit $q \to -m$ in **(66)** using $n=2$ we find

$$\sum_{\substack{l=1\\l\neq m}}\frac{\psi(1+l)}{(l-m)^2}=\psi(m+1)(\frac{1}{3}\pi^2-\psi'(m+1))-\psi^{(2)}(m+1)+\frac{\gamma}{m^2}. \quad \textbf{(69)}$$

Split the sum into its finite and infinite parts, apply the limit $q \to 0$ to (**C.6**) and compare the resulting infinite sums to obtain (**B.68**).

Variations of Euler (quadratic) sums are available by operating on **(47)** with $\lim\limits_{\mu \to q+1} \dfrac{\partial}{\partial q}$:

$$\sum_{l=0} \frac{\psi(q+l)\psi(q+l+1)}{(q+l)(1+l)} = \frac{1}{(q-1)}\Big(\frac{\pi^2}{6}[\psi(q)+\psi(q-1)]+\psi(q)\psi(q-1)[\psi(q)+\gamma] + 2\zeta(3)-\psi(q-1)\psi'(q)\Big) \qquad \textbf{(70)}$$

reducing – after considerable simplification - to a known result [**9**, Eq.(9)], [**16**, Eq.(8.92)] when $q \to 1$:

$$\sum_{l=0} \frac{[\psi(1+l)-\psi(1)]^2}{(1+l)^2} = \frac{11\,\pi^4}{360} \qquad \textbf{(71)}$$

Notice that the left-hand side of (**70**) can be rewritten as a transformation by a simple single step recursion, that is

$$\sum_{l=0} \frac{\psi(q+l)\psi(q+l+1)}{(q+l)(1+l)} = \sum_{l=0} \frac{\psi(q+l)^2}{(q+l)(1+l)} + \sum_{l=0} \frac{\psi(q+l)}{(q+l)^2(1+l)}\,, \qquad (\textbf{70}\text{a})$$

and so, in the limit $q \to 1$

$$\sum_{l=0} \frac{\psi(1+l)^2}{(1+l)^2} + \sum_{l=0} \frac{\psi(1+l)}{(1+l)^3} = \frac{\pi^2\gamma^2}{6} - 3\gamma\zeta(3) + \frac{\pi^4}{30} \qquad \textbf{(72)}$$

equivalent to Euler sums of weight (2,2) and (1,3), both of which are known (e.g. [**5**], [**9**]).

Alternatively, simple differentiation of **(47)** with respect to the variable *q* with (*p*>*q*+*1)* produces

$$\begin{aligned}\sum_{l=0} \frac{\psi(q+l)\psi(p+l)\Gamma(q+l)}{\Gamma(p+l)(1+l)} = \frac{\Gamma(q-1)}{\Gamma(p-1)}[&\psi(q-1)\{\psi'(p-q)-\psi'(p-1)\\ &+\psi(p-1)(\psi(p-1)-\psi(p-q))\}\\ &+\psi(p-1)\psi'(p-q)-\psi^{(2)}(p-q)]\,,\end{aligned} \qquad \textbf{(73)}$$

reducing to (**70**) when *p*=*q*+*1*. Similarly, operating on **(48)** with $\lim\limits_{p \to q+1} \dfrac{\partial}{\partial p}\dfrac{\partial}{\partial q}$ gives

$$\begin{aligned}\sum_{l=0} \frac{\psi(q+l)\psi(q+l+1)\psi(1+l)}{(q+l)(1+l)} =& \Big(\frac{1}{2}\psi(q)^4+\psi(q)^2\big(\frac{5\pi^2}{12}-\frac{3}{2}\psi'(q)-\frac{\gamma^2}{2}\big)\\ &+\psi(q)\big(\frac{1}{2}\psi^{(2)}(q)+4\zeta(3)\big)+\frac{11\pi^4}{180}-\frac{\pi^2}{6}\psi'(q)\Big)/(q-1)\\ &+\Big(-\frac{1}{2}\psi(q)^3+\psi(q)\big(\frac{3}{2}\psi'(q)-\frac{\pi^2}{4}+\frac{\gamma^2}{2}\big)-\frac{1}{2}\psi^{(2)}(q)-\zeta(3)\Big)/(q-1)^2\end{aligned} \qquad \textbf{(74)}$$

In the limit *q*=1, a straightforward result is

$$\sum_{l=0} \frac{\psi(2+l)\psi(1+l)^2}{(1+l)^2} = -\frac{1}{6}\gamma^3\,\pi^2 + 4\,\zeta(3)\,\gamma^2 - \frac{7}{72}\pi^4\,\gamma + \frac{1}{3}\zeta(3)\,\pi^2 + 6\,\zeta(5) \qquad \textbf{(75)}$$

Using this technique, along with **(3)**, variations in the choice of $\lim\limits_{p\to q+n}$ and partial fraction decomposition, all Euler sums up to cubic order become accessible. The use of computer algebra is recommended.

**Theorem 4.11**

| | |
|---|---|
| $$\sum_{l=0} \frac{\Gamma(c+l)\psi(1+l)^2}{\Gamma(p+l)} = -\frac{\Gamma(c)}{\Gamma(p)}\,{}_4F_3\left({1,1,1,p-c \atop 2,2,p}\,\middle\|\,1\right) + \frac{\Gamma(c)}{\Gamma(p-1)(p-1-c)}\Big[\frac{\pi^2}{6}+\gamma^2 + \psi(p-1-c)[2\gamma-\psi(p-1)-\psi(c)] + \psi'(p-1-c) + \psi(p-1-c)^2 - \psi(p-1)[2\gamma-\psi(c)]\Big]$$ | **(76)** |

and

| | |
|---|---|
| $$\sum_{l=0} \frac{\Gamma(c+l)\psi'(1+l)}{\Gamma(p+l)} = \frac{\Gamma(c)}{\Gamma(p)}\,{}_4F_3\left({1,1,1,p-c \atop 2,2,p}\,\middle\|\,1\right) - \frac{\Gamma(c)}{\Gamma(p-1)(p-1-c)}\big(\psi'(p-1) - [\psi(p-1)-\psi(p-1-c)][\psi(p-1)-\psi(c)]\big)$$ | **(77)** |

**Proof:**

Following the method of Theorem 4.5, operate on both sides of **(50)** and **(51)** with both $\lim\limits_{f\to 1}\lim\limits_{\beta\to f}\frac{\partial^2}{\partial\beta^2}$ and $\lim\limits_{\beta\to 1}\lim\limits_{f\to\beta}\frac{\partial^2}{\partial f^2}$. The resulting expressions contain a large number of terms, all of which can be evaluated (with computer algebra) using results given elsewhere in this paper, especially **(47)** and **(48)**. Add and subtract the resulting expressions. **§** For the case $p = c+1$, see **(C.20)**. Note: – the version 2014 statement of **(77)** contained a transcription error.

Seen from another viewpoint, subtracting (**76**) from **(77)** gives an alternative identification for a particular ${}_4F_3(1)$ specifically

| | |
|---|---|
| $${}_4F_3\left({1,1,1,c \atop 2,2,p}\,\middle\|\,1\right) = \frac{\Gamma(p)}{2\Gamma(p-c)}\sum_{l=0}\frac{\Gamma(p-c+l)}{\Gamma(p+l)}(\psi'(1+l)-\psi(1+l)^2) + \frac{(p-1)}{2(c-1)}(\psi_-\psi_+ + 2(\gamma-\psi(p-c))\psi_- + \gamma^2 + \pi^2/6 + \psi'_+), \quad p-c>-1 \text{ and } c>1,$$ | **(78)** |

where $\psi_\pm \equiv \psi(c-1) \pm \psi(p-1)$ after the substitution $c \to p-c$, the significance of which will be discussed elsewhere (to be published). See also **(16)** and **(B.43).**

Oct. 15, 2017.

## Acknowledgements

Although I had been thinking about updating this paper for some time, I am grateful to Prof. Lu Wei for motivating me to finally do so. I am also grateful to those (cited) authors who chose to make copies of their work freely and publicly accessible, and acknowledge the Researchgate computer algorithm which brought many of those papers to my attention. Some authors have chosen not to post freely available preprints or legally permitted offprints on public websites; citations to any such papers published after 2005 and only available behind a paywall, have been withheld. Please contact this author for further information.

## Appendix A: Sketch of Proof of (16)*

**Proof:** In **(65)**, with $\Re(c) > 0$, set n=0, giving

$$\sum_{l=0}[\frac{\psi(c+l)}{(l+1)}-\frac{\psi(l+1)}{(c+l)}]=\frac{\pi^2}{3}+\gamma\psi(c-1)-\psi'(c)+\psi(c)^2-\frac{\psi(c)}{(c-1)} \qquad \textbf{(A.1)}$$

Split the infinite sum into two parts: $l = 0\cdots k$, , and $l = k+1\cdots\infty$ , and shift the lower value of the summation index of the second (infinite) series to zero. The resulting infinite series, given by **(65)** with the replacement $c \to c+k+1$ and $n = k+1$, involves a $_4F_3$, so using **(A.1)** and **(21)**, solve for the $_4F_3$ and simplify. §

## Appendix B: Finite Sums

The following useful and related finite sums, based on the methods of this paper, many of which are needed [e.g. **31**]) or obtained elsewhere [**19**], are presented with minimal proof.

From [**21**, Eq.(3)]

$$\sum_{l=1}^{k}\frac{\psi(l+\frac{1}{2})}{(k-l+\frac{1}{2})}=2\sum_{l=1}^{k}\frac{\psi(l+\frac{1}{2})}{l}+\psi(\tfrac{1}{2})\left(H_k(\tfrac{1}{2})-2H_k(1)\right) \qquad \textbf{(B.1)}$$

From (**B.5**) (below), after some rearrangement and utilizing **(26)** obtain

$$\begin{aligned}\sum_{l=0}^{k}\frac{\psi(l+1)}{(k-l+\frac{1}{2})}=&-\sum_{l=0}^{k}\frac{\psi(l+\frac{1}{2})}{l+1}+\tfrac{1}{2}\left(H_{k+1}(\tfrac{1}{2})^2+H'_{k+1}(\tfrac{1}{2})\right)\\&+H_{k+1}(1)\left(\psi(k+\tfrac{3}{2})+2\right)-(\gamma+2)H_{k+1}(\tfrac{1}{2})\end{aligned} \qquad \textbf{(B.2)}$$

and

$$\sum_{l=0}^{k}\frac{\psi(l+1)}{(k-l+1)}=\gamma\psi(k+2)+\psi'(k+2)+\psi(k+2)^2-\frac{\pi^2}{6} \qquad \textbf{(B.3)}$$

which, after reversal, or with reference to **(21)**, can also be written

$$\sum_{l=1}^{k}\frac{\psi(k-l+1)}{l}=\gamma\psi(k+1)+\psi'(k+1)+\psi(k+1)^2-\frac{\pi^2}{6}\,.$$

Also, see the case $\beta \to n$ in (**B.53**) below. Equivalently,

* A proof of **(16)** based on [**26**, Eq. 7.4.4(38)] might be attainable. However, that result fails several numerical tests, and appears to be incorrect.

$$\sum_{l=0}^{k}\frac{1}{(k-l+\frac{1}{2})}\sum_{m=0}^{l}\frac{1}{(1+m)}=-\sum_{l=0}^{k}\frac{1}{(l+1)}\sum_{m=0}^{l}\frac{1}{(m+\frac{1}{2})}$$
$$+\tfrac{1}{2}\left(H_{k+2}(\tfrac{1}{2})^2+H'_{k+2}(\tfrac{1}{2})\right)+H_{k+1}(1)H_{k+2}(\tfrac{1}{2}) \quad \textbf{(B.4)}$$

and from [**21**] (also see **(B.13)**)

$$\sum_{l=0}^{k}\frac{1}{k-l+c}\sum_{m=0}^{l}\frac{1}{m+c}=2\sum_{l=0}^{k}\frac{1}{l+2c}\sum_{m=0}^{l}\frac{1}{m+c}. \quad \textbf{(B.5)}$$

Follow the same steps as in **(22)** and operate with $\frac{\partial^n}{\partial b^n}$ to generalize **(21)**:

$$\sum_{l=0}^{k}\frac{\psi^{(n)}(b+l)}{(c+l)^s}=\psi^{(n)}(b)\zeta(s,c)-\psi^{(n)}(b+k+1)\zeta(s,c+k+1)$$
$$+(-1)^n n!\sum_{l=0}^{k}\frac{\zeta(s,c+l+1)}{(b+l)^{n+1}}, \quad \textbf{(B.6)}$$

where $s$ is arbitrary, except $s\neq 1$. In the case $n=1, s=1$, differentiate **(21)** once with respect to the parameter $b$ to find

$$\sum_{l=0}^{k}\frac{\psi'(b+l)}{c+l}=\sum_{l=0}^{k}\frac{\psi(c+l)}{(b+l)^2}+\psi(c+k+1)\psi'(b+k+1)+\frac{\psi'(b+k+1)-\psi'(b)}{(b-c)}$$
$$-\psi'(b)\psi(c)+\frac{\psi(c+k+1)-\psi(b+k+1)+\psi(b)-\psi(c)}{(b-c)^2}. \quad \textbf{(B.7)}$$

In the limit $c\to b$, we obtain

$$\sum_{l=0}^{k}\frac{\psi'(b+l)}{b+l}=\sum_{l=0}^{k}\frac{\psi(b+l)}{(b+l)^2}+\frac{1}{2}\ \psi^{(2)}(b+k+1)+\psi(b+1+k)\psi'(b+1+k)$$
$$-\psi'(b)\psi(b)-\frac{1}{2}\psi^{(2)}(b)\ . \quad \textbf{(B.8)}$$

Further, in **(B.6)** set $b=c,\ \ s=n+1$ and solve, thereby producing a variant of **(26)**:

$$\sum_{l=0}^{k}\frac{\psi^{(n)}(b+l)}{(b+l)^{n+1}}=\frac{(-1)^n}{2}\Big\{\frac{n!}{(2n+1)!}[\psi^{(2n+1)}(b+1+k)-\psi^{(2n+1)}(b)]$$
$$+\frac{1}{n!}[\psi^{(n)}(b+1+k)^2-\psi^{(n)}(b)^2]\Big\}\ , \quad \textbf{(B.9)}$$

Noting that **(21)** is valid for all $c$, let c $\to -k-1$, evaluate that limit, compare the coefficients of the first order expansion terms and replace c $\to b$ to find

$$\sum_{l=0}^{k}\frac{\psi(b+l)}{(k-l+1)^2}=\frac{\pi^2}{6}\psi(b+k+1)-\psi'(k+2)\psi(b)-\sum_{l=0}^{k}\frac{\psi'(k-l+1)}{b+l}. \qquad \textbf{(B.10)}$$

See also (**B.7**). An equivalent form of **(B.10)** can be found by reversing each of the sums, and replacing $b \to b-k$ to yield

$$\sum_{l=0}^{k}\frac{\psi(b-l)}{(l+1)^2}=\frac{\pi^2}{6}\psi(b+1)-\psi'(k+2)\psi(b-k)-\sum_{l=0}^{k}\frac{\psi'(l+1)}{b-l}. \qquad \textbf{(B.11)}$$

The case **(B.9)** can also be used to find more general results. For example, after partial fraction decomposition and some rearrangement, the case $n$=1 using (**B.8**) and **(B.9)** gives

$$\begin{aligned}\sum_{l=0}^{k-1}\frac{\psi'(1+l)}{(1+l)^2(2+l)^2}&=-\psi'(k)^2-\psi'(k)\frac{(4k^4+6k^3-k^2-4k-2)}{k^2(k+1)^2}\\&\quad-\frac{\psi^{(3)}(k)}{6}-\frac{(5k+1)(k-1)}{k^2}+\frac{7\pi^4}{180}+\frac{\pi^2}{6}\end{aligned}. \qquad \textbf{(B.12)}$$

Also, simplifying [**21**, Eqs. (3) and (5)] yields

$$\sum_{l=0}^{k}\frac{\psi(l+1+c)}{(2c+l)}=\frac{1}{2}\sum_{l=0}^{k}\frac{\psi(l+1+c)}{(k-l+c)}+\left(\psi(2c+k+1)-\psi(2c)-\tfrac{1}{2}\psi(k+c+1)\right)\psi(c)+\tfrac{1}{2}\psi(c)^2. \qquad \textbf{(B.13)}$$

In the case that $c=-k$, an extensive series of limits utilizing **(26)** and (**B.44**) among others, eventually reduces (**B.13**) to (**B.10**). Note that (**B.16**) below is a special case of (**B.13**), both in the limit $c \to 0$. Also see **(B.66)**.

Several new results can be obtained by comparing **(16)** with **[31**, Eq. (A.24) (generalizing $n \to c$ )**]** to obtain

$$\sum_{l=1}^{m}\frac{\psi(c-l+1)}{l}+\sum_{l=1}^{m}\frac{\psi(c-m+l)}{l}=(\gamma+\psi(m+1))\psi_{+}+\tfrac{1}{2}(\psi_{-}^{\;2}-\psi_{-}') \qquad \textbf{(B.14)}$$

where $\psi_{\pm}=\psi(c-m)\pm\psi(c+1)$. Further, reversing either term in the left-hand side of (**B.14**) identifies the related sums

$$\sum_{l=1}^{m}\frac{\psi(c-l+1)}{l}+\sum_{l=1}^{m}\frac{\psi(c-m+l)}{l}=(m+1)\sum_{l=1}^{m}\frac{\psi(c-m+l)}{l(m+1-l)}=(m+1)\sum_{l=1}^{m}\frac{\psi(c+1-l)}{l(m+1-l)} \qquad \textbf{(B.15)}$$

and a similar result exists for (**B.16**) below, if the summations are combined. After applying **(19)** to (**B.14**) we find another variation

$$\sum_{l=1}^{m}\left(\frac{\psi(l)}{c-l+1}-\frac{\psi(l)}{c-m+l}\right)=\frac{1}{2}\left(\psi(c-m+1)-\psi(c+1)\right)^2-\frac{1}{2}\left(\psi'(c-m+1)-\psi'(c+1)\right) \qquad \textbf{(B.16)}$$

The special case $c = m/2$ in (**B.16**) requires a careful evaluation of limits. To first order, the case $m = 2j - 1$, yields

$$\sum_{l=1}^{2j-1}\frac{\psi(l-j+\frac{1}{2})}{l}=\psi(\tfrac{1}{2}+j)[\gamma+\psi(2j)]+\tfrac{1}{2}\psi'(\tfrac{1}{2}+j)-\frac{\pi^2}{4}\ , \qquad \textbf{(B.17)}$$

which can also be written (by splitting the sum and reversing one of the resultants)

$$\sum_{l=0}^{j-1}\frac{\psi(l+1/2)}{l^2-j^2}=-\frac{\psi(\frac{1}{2}+j)}{2j}[\gamma+\psi(2j)]-\frac{\psi'(\frac{1}{2}+j)}{4j}+\frac{\pi^2}{8j}-\frac{\psi(\frac{1}{2})}{2j^2}\ . \qquad \textbf{(B.18)}$$

Again, to first order, in the case that $m = 2j$, a very careful evaluation of limits yields

$$\sum_{l=1}^{j}\frac{\psi(l+j+1)}{l}=\psi(1+j)\gamma+\psi(j+1)^2-\tfrac{1}{2}\psi'(j)+\frac{\pi^2}{12}+\frac{1}{2j^2}\ , \qquad \textbf{(B.19)}$$

which, by virtue of **(23)**, can also be written

$$\sum_{l=0}^{j}\frac{\psi(l+1)}{l+j+1}=\psi(1+j)\psi(2+2j)-\psi(j+1)^2+\tfrac{1}{2}\psi'(j+1)-\frac{\pi^2}{12}\ . \qquad \textbf{(B.19a)}$$

By incrementing the argument of the digamma function on the left-hand side and identifying most of the resulting sums, an extension of **(B.19)** is found

$$\begin{aligned}\sum_{l=1}^{j}\frac{\psi(l+n)}{l}&=-\sum_{l=0}^{n-j-2}\frac{\psi(j+n-l)}{n-1-l}+\tfrac{1}{2}(\psi(1+j)+\psi(n))^2-\psi(j+1)^2\\&\quad-\tfrac{1}{2}\psi'(n)+\gamma\psi(n)+\frac{\pi^2}{12}\quad\text{if}\quad n\geq j+2.\end{aligned} \qquad \textbf{(B.20)}$$

In the case that $n = 1$, **(B.20)** reduces to the known result

$$\sum_{l=1}^{j}\frac{\psi(l+1)-\psi(1)}{l}=\tfrac{1}{2}(\psi(1+j)+\gamma)^2-\tfrac{1}{2}\psi'(j+1)+\frac{\pi^2}{12}=-\int_0^1\frac{1-t^j}{1-t}\log(1-t)dt, \qquad \textbf{(B.21)}$$

the integral arising from Hansen [**12**, Eq. (55.2.1)]. For the case $n = 0$, see (**B.66**). Note that, if the upper limit becomes negative in **(B.20),**

$$\begin{aligned}\sum_{l=0}^{n-j-2}\frac{\psi(j+n-l)}{n-1-l}&=\sum_{l=n-j-1}^{-1}\frac{\psi(j+n-l)}{l+1-n}\quad\text{if } n\leq j\\&=0\qquad\qquad\text{if } n=j+1\ .\end{aligned} \qquad \textbf{(B.22)}$$

Considering the second order terms of the case $c = m/2$ in (**B.16**) leads to

$$\sum_{l=1}^{j}\frac{\psi'(j-l+1)}{l}+\sum_{l=1}^{j}\frac{\psi'(l)}{l}=2\psi'(1+j)(\gamma+\psi(j+1))+2\zeta(3)+\psi^{(2)}(1+j). \qquad \textbf{(B.23)}$$

Or, if $c = j - 1$ in (**B.16**), the second order terms lead to

$$\sum_{l=0}^{j-2}\frac{\psi(l+1)}{(1+l-j)^3}+\sum_{l=0}^{j-2}\frac{\psi(l+j+1)}{(l+1)^3}=\frac{\pi^2}{3}\psi'(j)-\psi'(j)^2-\frac{5}{12}\psi^{(3)}(j) \qquad (\textbf{B.23}\text{a})$$

For other variations reverse the sums, or apply (**B.44**) below.

From the literature **[37, Eq. (1.6)**, corrected for the wrong sign in the first term of the right-hand side] by considering even and odd terms, and after applying (**B.67** – see below) we obtain

$$\begin{aligned}\sum_{l=1}^{n-1}\frac{\psi(\frac{1}{2}l)}{l-n}=(-)^n\Big(&\sum_{k=1}^{\lfloor\frac{n+1}{2}\rfloor}\frac{\psi(l)}{l-\frac{1}{2}}+(\ln(2)-\psi(n)+\tfrac{1}{2}\psi(\tfrac{1}{2}n))\psi(\tfrac{1}{2}n+1)+\tfrac{1}{2}\psi(\tfrac{1}{2})^2-\frac{\pi^2}{12}\Big)\\&+\psi(\tfrac{1}{2})(2\psi(\tfrac{1}{2}n)+\ln(2)-\psi(n))+\ln(2)(2\psi(\tfrac{1}{2}n)-\tfrac{2}{n})-\psi(\tfrac{1}{2}n)\psi(n+1)\\&-\tfrac{1}{2}\psi'(\tfrac{1}{2}n)+\frac{\pi^2}{6}+\tfrac{2}{n}\psi(n),\end{aligned} \qquad \textbf{(B.24)}$$

and, from **[37, Eq.(1.8)]** using the same method, we find

$$\begin{aligned}\sum_{l=1}^{n-1}\frac{\psi(l+\frac{1}{2})}{l-n}=&\sum_{l=0}^{n-1}\frac{\psi(l+\frac{1}{2})}{l+1}+\tfrac{1}{2}\psi(\tfrac{1}{2})^2-(3\psi(n+1)-4\psi(2n+1)+2\ln(2)-\tfrac{1}{n})\psi(\tfrac{1}{2})\\&-6\ln^2(2)-(2\psi(n+1)-8\psi(2n+1))\ln(2)+\tfrac{1}{2}\psi(n+1)^2-4\psi(n+1)-2\psi(2n+1)^2\\&+4\psi(2n+1)-\tfrac{1}{2}\psi'(\tfrac{1}{2}+n)+\frac{\pi^2}{4}.\end{aligned} \qquad \textbf{(B.25)}$$

Operating on (**28**) with $\dfrac{\partial}{\partial a}$, and evaluating the limit $a \to b+1$ leads to

$$\begin{aligned}\sum_{l=0}^{k}\frac{\psi(b+l)\psi(b+l+1)}{(b+l)}&=\sum_{l=0}^{k}\Big\{\frac{\psi(b+l)^2}{b+l}+\frac{\psi(b+l)}{(b+l)^2}\Big\}\\&=\tfrac{1}{3}(\psi(b+k)^3-\psi(b)^3)+\frac{\psi(b+k)^2}{(b+k)}+\frac{\psi(b+k)}{(b+k)^2}+\tfrac{1}{6}(\psi^{(2)}(b)-\psi^{(2)}(b+k))\,.\end{aligned} \qquad \textbf{(B.26)}$$

Use **(B.9)** with $n$=0, $s$=2 to obtain a variant of the above – see (**B.8**). For the case $a \to b$ see (**B.31**).

As has been noted, operating on well-known hypergeometric sums is a technique extensively utilized in the literature [e.g. **17, 8**]. Typically, such calculations are based on variants of Gauss' evaluation of any ${}_2F_1(1)$ and the well-known Dixon/Whipple/Watson evaluation(s) of a particular set of ${}_3F_2(1)$ (e.g. [**29**], [**33**] ,[**38**]). Here, it is noted that other known ${}_3F_2(1)$ evaluations lend themselves well to this technique. Consider the tabulated Entry 36 (see Milgram, [**23**]), specifically ${}_3F_2$(1,1-n, a;b,c;1) satisfying

$$\sum_{l=0}^{n-1}\frac{(-1)^l\Gamma(a+l)}{\Gamma(n-l)\Gamma(b+l)\Gamma(c+l)}=\frac{\Gamma(1+a-c)\Gamma(a)\sin(\pi c)}{\pi\ \Gamma(b-1+n)\Gamma(b-1)}\sum_{l=0}^{n-1}\frac{\Gamma(2+l-c-n)\Gamma(b-1+l)}{\Gamma(l+1)\Gamma(a-c+2-n+l)} \qquad \textbf{(B.27)}$$

Apply the ordered operator $\lim\limits_{c\to 1-n}\left(\lim\limits_{a\to b+1}\frac{\partial}{\partial c}\frac{\partial}{\partial a}\right)$ to both sides of (**B.27**) and, after considerable calculation involving the sums listed herein, eventually arrive at

$$\begin{aligned}\sum_{l=0}^{n-1}\frac{\psi(a+l)}{(a+l)}(\psi(l+1)+\psi(n-l)) &= -\sum_{l=0}^{n-1}\frac{\psi(l+1)}{(a+l)^2}+\frac{1}{6}\psi(a)^3\\ &+(\psi(a)+\psi(a+n))\sum_{l=0}^{n-1}\frac{\psi(l+1)}{(a+l)}+\frac{1}{3}\psi(a+n)^3+(\psi'(a+n)-\psi'(a))\psi(a+n)\\ &-\frac{1}{2}\psi(a)\psi(a+n)^2+\frac{1}{2}(\psi'(a)-\psi'(a+n))\psi(a)+\frac{1}{3}(\psi^{(2)}(a+n)-\psi^{(2)}(a))\end{aligned} \qquad \textbf{(B.28)}$$

In the case that $a=(1-n)/2$, (**B.28**) reduces to an equivalent, but simpler form of **(47)** when $k=2j$, that being

$$\sum_{l=1}^{2j}\frac{\psi(l)}{\frac{1}{2}+j-l}=\frac{1}{2}\psi'(\tfrac{1}{2}+j)-\pi^2/4\ . \qquad \textbf{(B.29)}$$

If $a=1-j$, (**B.28**) reduces to **(18).** Alternatively, apply the ordered operator $\lim\limits_{c\to 1-n}\left(\lim\limits_{b\to a}\frac{\partial}{\partial a}\right)$ to (**B.27**) and again, after a long series of simplifications involving many of the sums listed here, as well as the following result, (equivalent to [**32**, Section 2, Theorem 1 with $l=1$ ], and generalizing [**31**, Eq.(A.7)] by the replacement $a\to k$ )

$$\sum_{l=0}^{n-1}\psi'(a+l)=\psi(a+n)-\psi(a)+(1-a)\psi'(a)+(a+n-1)\psi'(a+n), \qquad \textbf{(B.30)}$$

find

$$\begin{aligned}\sum_{l=0}^{n-1}\psi(a+l)^2 &= (a+n-1)\psi(a+n)^2+(1-2(a+n))\psi(a+n)\\ &+(1-a)\psi(a)^2+(2a-1)\psi(a)+2n\ .\end{aligned} \qquad \textbf{(B.31)}$$

For related results of the above form see Wei [**31**, Eq. (A.4)] and Wei et. al. [**32**], and for those forms involving alternating series with binomial coefficients, but lacking the independent parameter "$a$", see Choi [**4**].

Similarly, consider the evaluated Entry 25 of Milgram, [**23**], that is ${}_3F_2$(a,1,b ; m+1,c;1) (originally **(12)**). Operate on both sides with $\lim\limits_{c\to 1+a}\frac{\partial}{\partial a}\left(\lim\limits_{b\to 1-m}\frac{\partial}{\partial b}\right)$ to eventually obtain

$$\sum_{l=0}^{m-1}\frac{\Gamma(a-l)\psi(1+l)}{\Gamma(m-l)(l+1)}=\frac{\Gamma(a+1)}{\Gamma(m+1)}\Big(\gamma(\psi(a-m+1)-\psi(a+1))$$
$$+\frac{1}{2}(\psi(a-m+1)-\psi(a+1))^2-\frac{1}{2}(\psi'(a-m+1)-\psi'(a+1))\Big). \quad \textbf{(B.32)}$$

See also **(B.27)** and **(B.50)** as well as Hansen [**12**, section 55.4]. Note that the operator $\left(\lim_{a\to m}\frac{\partial}{\partial a}\right)$ applied to (**B.32**) yields (**B.58**) below.

A related sum can be found based upon a variation of **(16**), as well as **(18**), that being

$$\sum_{l=0}^{m-1}\frac{\Gamma(a-l)\psi(1+l)\psi(2+l)}{\Gamma(m-l)(l+1)}=\frac{\Gamma(a+1)}{\Gamma(m+1)}\Big(\frac{1}{6}\Psi_{(0)}^3+\frac{g}{2}(\Psi_{(0)}^2-\Psi_{(1)})$$
$$+(\gamma+\psi(a+1)+\Psi_{(0)})\sum_{l=0}^{m-1}\frac{\psi(1+l)}{a-m+l+1}-\sum_{l=0}^{m-1}\frac{\psi(a-l)\psi(1+l)}{a-l} \quad \textbf{(B.33)}$$
$$-\frac{1}{3}\Psi_{(2)}+(-\gamma^2+\gamma\psi(m+1)+\frac{1}{2}\Psi_{(1)})\Psi_{(0)}\Big)$$

where (also see **(9)**)

$$\Psi_{(n)}=\psi^{(n)}(a-m+1)-\psi^{(n)}(1+a) \text{ and } g=-\gamma+\psi(m+1)+\psi(a+1)\,.$$

A related result follows:

$$\frac{(\psi(m+1)-\gamma)\Gamma(m+1)}{\Gamma(a+1)}\sum_{l=0}^{m-1}\frac{\Gamma(a-l)\psi(1+l)}{\Gamma(m-l)(l+1)}=(\psi(a+1)-\gamma)\sum_{l=0}^{m-1}\frac{\psi(1+l)}{a-m+l+1}$$
$$-\sum_{l=0}^{m-1}\frac{\psi(a-m+l+1)\psi(1+l)}{a-m+l+1}-\sum_{l=0}^{m-1}\frac{\psi(1+l)}{(a-m+l+1)^2} \quad \textbf{(B.34)}$$
$$+\frac{\Gamma(m+1)}{\Gamma(a+1)}\sum_{l=0}^{m-1}\frac{\Gamma(a-l)\psi(1+l)\psi(2+l)}{\Gamma(m-l)(l+1)}$$

Eliminating terms between **(B.33**) and **(B.34**) is equivalent to **(B.32**). Similar results may be obtained from the well-known, simple hypergeometric identities (see also **(B.61**))

$$\sum_{l=0}^{n-1}\frac{\Gamma(l-a)}{(l-n)\Gamma(l+1)}=\frac{\Gamma(n-a)}{\Gamma(1+n)}(\psi(-a)-\psi(n-a)) \quad \textbf{(B.35)}$$

and

$$\sum_{l=0}^{n-1}\frac{(-1)^l}{\Gamma(n-l)\Gamma(l-a+1)}=\frac{a}{\Gamma(n)\Gamma(1-a)(1-n+a)}. \quad \textbf{(B.36)}$$

After differentiating once with respect to *a*, the first yields

$$\sum_{l=0}^{n-1}\frac{\Gamma(l-a)\psi(l-a)}{(l-n)\Gamma(l+1)}=\frac{\Gamma(n-a)}{\Gamma(1+n)}\Big(\psi(-a)\psi(n-a)-\psi(n-a)^2-\psi'(n-a)+\psi'(-a)\Big)\,, \quad \textbf{(B.37)}$$

and the second yields the known result

$$\sum_{l=0}^{n-1}\frac{(-1)^l\psi(l-a+1)}{\Gamma(n-l)\Gamma(l-a+1)}=\frac{a\psi(1-a)(1-n+a)+1-n}{(1-n+a)^2\,\Gamma(n)\Gamma(1-a)}. \quad \textbf{(B.38)}$$

Comparing the zero^th order terms in the limiting case $a\to n-1$ gives the known result

$$\sum_{l=0}^{n-1}\frac{\psi(l+1)}{l+1}=\sum_{l=0}^{n-1}\frac{\psi(n-l)}{n-l}=\tfrac{1}{2}\psi(n)^2+\psi(n)/n+\tfrac{1}{2}\psi'(n)-\pi^2/12-\tfrac{1}{2}\gamma^2. \quad \textbf{(B.39)}$$

Additionally, the first order terms in this limiting case give

$$\sum_{l=0}^{n-2}\frac{\psi(n-l-1)^2}{1-n+l}+\sum_{l=0}^{n-2}\frac{\psi'(n-l-1)}{1-n+l}=-\psi(n)\psi'(n-1)-\frac{\psi(n-1)^2}{(n-1)}$$
$$-\frac{1}{3}(\psi(n-1)^3+\psi^{(2)}(n-1))-\frac{\pi^2\gamma}{6}-\frac{1}{3}\gamma^3-\frac{2\zeta(3)}{3}. \quad \textbf{(B.40)}$$

With respect to the former case, see (**23)** and **(B.9**); otherwise, see **(B.44**), **(52)** and **(B.59);** regarding either case, reverse the order as illustrated throughout.

Alternatively, consider the aforementioned Entry 25 under the ordered substitution $n=1,\ b\to 1-n, c\to 1-c$**.** Applying the ordered operator $\left(\lim_{c\to a}\frac{\partial}{\partial c}\frac{\partial}{\partial a}\right)$ yields

$$\sum_{l=0}^{n-1}\frac{(-1)^l\psi(a+l)^2}{\Gamma(n-l)\Gamma(l+2)}=\frac{\Gamma(a-1)}{n\Gamma(a-1+n)}(\psi(a-1)+\psi(a-1+n)-2\gamma-2\psi(n))+\frac{\psi(a-1)^2}{\Gamma(n+1)} \quad \textbf{(B.41)}$$

Or, differentiate twice with respect to the parameter "a", followed by $a\to 1$ to obtain

$$\sum_{l=0}^{n-1}\frac{\psi(l+1)^2\,\Gamma(-c-l)}{(l+1)\Gamma(n-l)}=-\sum_{l=0}^{n-1}\frac{\psi'(l+1)\Gamma(-c-l)}{(1+l)\Gamma(n-l)}$$
$$+\frac{\Gamma(1-c)}{\Gamma(n+1)}\left(\tfrac{1}{3}(\psi_C{}^3+\psi_C{}^{(2)})-\gamma\psi_C{}^2+\psi_C[\gamma^2+\tfrac{\pi^2}{6}-\psi'_C]+\gamma\psi'_C\right) \quad \textbf{(B.42)}$$

where

$$\psi_C=\psi(c)-\psi(c+n)\,.$$

In the limit $c\to -n$, (**B.42**) reduces to a variant of (**B.58**), but also see (**B.47**).

Using Entry 36 of Milgram, [**23**], that is $_3F_2(1,1\text{-m, b ; c, 1-}\beta\text{ ;1})$, apply the ordered operator $\left(\lim_{c\to 1}\lim_{b\to 1}\frac{\partial}{\partial c}\frac{\partial}{\partial b}\right)$ to obtain

$$m_\beta\,\Gamma(-\beta)\Gamma(m)\sum_{l=0}^{m-1}\frac{(-1)^l\psi(1+l)^2}{\Gamma(m-l)\Gamma(l-\beta+1)}=-\sum_{l=0}^{m-1}\frac{\psi(m_\beta+1+l)}{(l+1)}-\gamma^2-\frac{\gamma}{m_\beta}$$
$$-\frac{1}{2}(\psi'(1+\beta)-\psi'(m_\beta))+(\psi(m)-\gamma-\frac{2}{m_\beta})\psi(1+m_\beta) \qquad \textbf{(B.43)}$$
$$+(2\gamma+\frac{2}{m_\beta}+\frac{1}{m})\psi(1+\beta)-\frac{\psi(m)}{m_\beta}-\frac{1}{2}(\psi(m_\beta)-\psi(1+\beta))^2$$

where $m_\beta=1-m+\beta$. In the case $m_\beta\to 0$, equating the coefficients of the zero[th] and first order terms respectively gives

$$\sum_{l=0}^{m-1}\frac{\psi(1+l)^2}{l-m}=-\sum_{l=0}^{m}\frac{\psi'(1+l)}{(l+1)}+\psi'(m+1)(\gamma+\frac{1}{m+1}-\psi(m+1))-\psi(m+1)^3-\psi(m+1)^2\gamma+\frac{\pi^2}{3}\psi(m+1) \qquad \textbf{(B.44)}$$

and

$$\sum_{l=0}^{m-1}\frac{\psi(l+1)^2\psi(m-l)}{m-l}=\psi(m+1)\sum_{l=0}^{m}\frac{\psi'(l+1)}{l+1}+\frac{1}{2}\sum_{l=0}^{m}\frac{\psi^{(2)}(l+1)}{l+1}+\frac{\pi^4}{120}$$
$$+\psi(m+1)^2(3\psi'(m+1)-5\pi^2/6-\gamma^2)/2+\psi(m+1)^4/2+\psi^{(3)}(m+1)/12$$
$$-\psi(m+1)\Big((\gamma+1/(m+1))\psi'(m+1)-\psi^{(2)}(m+1)/2-2\zeta(3)\Big) \qquad \textbf{(B.45)}$$
$$+\psi'(m+1)^2/2-\pi^2\psi'(m+1)/6-\big(\gamma+1/(m+1)\big)\psi^{(2)}(m+1)/2\ .$$

By differentiating (**B.43**) once with respect to $\beta$ , a lengthy, but general, result for the sum

$$\sum_{l=0}^{m}\frac{\psi(l+1)^2\psi(l-\beta+1)}{\Gamma(m-l+1)\Gamma(l-\beta+1)} \qquad (\textbf{B.43}\text{a})$$

can be found. Because of its complexity, that full result (whose simplification, if any exists, eludes me) is presented in the form of computerized output – see Appendix **D**. In the case $\beta=m+1$, we find

$$\sum_{l=0}^{m}\psi(l+1)^2\psi(m-l+1)=(m+1)\sum_{l=0}^{m}\frac{\psi'(l+1)}{(l+1)}+\frac{1}{3}(m+1/2-(m+1)\psi(m+1))\pi^2$$
$$-(m+1)\big(\gamma+2-\psi(m+1))\psi'(m+1)\big)+(m+1)\psi(m+1)^3-(3m+1)\psi(m+1)^2 \qquad \textbf{(B.46)}$$
$$+(1-\gamma+6m)\,\psi(m+1)+\gamma-6m\ ,$$

and for the case $\beta = m$ see (**B.45**).

Obvious variations of the above can be obtained by reversing the sum on the left-hand side or applying any combination of (**B.8**), **(B.10)** and/or **B.58**); note also the parallels between (**B.44**) and **(B.16)** taking **(21)** into account or (**B.28**) when *a=1*. Also, comparing (**B.44**) with (**B.42**) in the case that $c\to -n$ leads to

$$\sum_{l=0}^{n-2}\frac{\psi(1+l)^2}{l-n+1}-\sum_{l=0}^{n-2}\frac{\psi(1+l)^2}{(l+1)}=-\psi(n)^2\gamma+\frac{1}{6}(-2\psi(n)+\gamma)(-\pi^2+6\psi'(n))$$
$$-\frac{\gamma^3}{3}-\frac{2}{3}\zeta(3)-\frac{1}{3}\psi^{(2)}(n)-\frac{4}{3}\psi(n)^3\ . \qquad \textbf{(B.47)}$$

In the limit $n\to\infty$ the difference of the sums in (**B.47**) diverges as $-\frac{4}{3}\ln^3(n)$. Further, **(21)** with $b=1$ applied to (**B.47**) leads to

$$\sum_{l=0}^{n-2}\frac{\psi(1+l)^2}{l-n+1}+\sum_{l=0}^{n-2}\frac{\psi(1+l)}{(l+1)^2}=-\psi(n)^2\gamma+\frac{1}{6}(-2\psi(n)+\gamma)(-\pi^2+6\psi'(n))$$
$$-\zeta(3)-\frac{1}{2}\psi^{(2)}(n)-\psi(n)^3\ ,\quad n\ge 2 \qquad \textbf{(B.48)}$$

Compare with (**B.40**). Further variations can be identified by splitting any of the above the sums according to whether $n$ is even or odd. Specifically, in the case of (**B.47**), for $n>0$ and $n\to 2n$ we have

$$\sum_{l=0}^{2n-2}\frac{\psi(1+l)^2}{l-2n+1}-\sum_{l=0}^{2n-2}\frac{\psi(1+l)^2}{(l+1)}=2n\Big(\sum_{l=0}^{n-1}\frac{\psi(1+l)^2}{(l-2n+1)(l+1)}+\sum_{l=0}^{n-2}\frac{\psi(l+n+1)^2}{(l+1)^2-n^2}\Big), \qquad \textbf{(B.49)}$$

and, in the case that $n\to 2n+1$ we find

$$\sum_{l=0}^{2n-1}\frac{\psi(1+l)^2}{l-2n}-\sum_{l=0}^{2n-1}\frac{\psi(1+l)^2}{(l+1)}=(2n+1)\Big(\sum_{l=0}^{n}\frac{\psi(1+l)^2}{(l-2n)(l+1)}+\sum_{l=0}^{n-2}\frac{\psi(l+n+2)^2}{(l+3/2)^2-(n+1/2)^2}\Big)\ . \qquad \textbf{(B.50)}$$

Another method that can be used to obtain useful sums is to make use of known hypergeometric (Thomae) transformations of any ${}_3F_2(1)$. Consider Entry 25 -see **(B.42)**- and Entry 26 – ${}_3F_2$(a, n+1,b ; n+2,c;1) - (both listed in [**23**]) which can be transformed into one another by careful use of limits. Let $c=\alpha, b=\beta, n\to n+1$ in both entries, and evaluate the limit $a\to n+1$ in Entry 25 and the limit $a\to 1$ in Entry 26. The left-hand sides of both entries are now equal, so equating the right-had sides of each and replacing $\alpha\to-\alpha+\beta$ gives a useful identity with two free parameters:

$$\sum_{l=0}^{n}\frac{(-1)^l\psi(l+1)\Gamma(1+\alpha+l)}{\Gamma(l+1)\Gamma(n-l+1)\Gamma(\beta+l-n)}=$$
$$\frac{-\Gamma(\alpha+1)}{\Gamma(\beta)\Gamma(n+1)\Gamma(2+\alpha-\beta)}\Big(\frac{\pi}{\sin(\pi\beta)\Gamma(\beta-n-1)}\sum_{l=0}^{n-1}\frac{\Gamma(2+l+\alpha-\beta)}{(l-n)\Gamma(2+l-\beta)}$$
$$-(-1)^n\Gamma(2+n+\alpha-\beta)\big(\psi(\beta-\alpha-n-1)-\psi(\beta-\alpha-1)+\psi(n+1)\big)\Big). \qquad \textbf{(B.51)}$$

In the case $\alpha=0$ , $\beta\to n+1-\beta$ , (**B.51**) reduces to

$$\sum_{l=0}^{n}\frac{(-1)^l\psi(l+1)}{\Gamma(n-l+1)\Gamma(-\beta+l+1)}=\frac{\psi(-\beta)-\psi(-\beta+n)-\gamma}{(-\beta+n)\Gamma(-\beta)\Gamma(n+1)}. \quad \textbf{(B.52)}$$

Differentiate (**B.52**) once with respect to β to obtain

$$\sum_{l=0}^{n}\frac{(-1)^l\psi(l+1)\psi(l-\beta+1)}{\Gamma(n-l+1)\Gamma(-\beta+l+1)}=\frac{-\big(\gamma-\psi(-\beta)+\psi(n-\beta)\big)\big(\psi(-\beta)+1/(n-\beta)\big)-\psi'(-\beta)+\psi'(-\beta+n)}{(-\beta+n)\Gamma(-\beta)\Gamma(n+1)} \quad \textbf{(B.53)}$$

and, in the limit $\beta\to n-1$, comparison of the three lowest order coefficients yields

$$\sum_{l=0}^{n}\psi(l+1)=(\psi(n+2)-1)(n+1), \quad \textbf{(B.54)}$$

$$\sum_{l=0}^{n}\psi(l+1)\psi(n-l+1)=(n+1)(\psi(n+2)^2+\psi'(n+2)-2\psi(n+2)+2-\pi^2/6) \quad \textbf{(B.55)}$$

and

$$\sum_{l=0}^{n}\psi(l+1)\psi(n-l+1)^2+\sum_{l=0}^{n}\psi(l+1)\psi'(n-l+1)=$$
$$\frac{(n+1)}{3}(3\psi(n+2)^3-9\psi(n+2)^2+\psi(2+n)\big(18+9\psi'(2+n)-\pi^2\big)$$
$$-9\psi'(2+n)+3\psi^{(2)}(2+n)+6\zeta(3)+\pi^2-18). \quad \textbf{(B.56)}$$

(**B.54**) is a well-known result [see Hansen **12**, Eq. (55.6.1)], while (**B.55**) has the interesting property that the left-hand side is invariant under reversal, so, from **(67a)** we find

$$\sum_{l=0}^{n}l\,\psi(l+1)\psi(n-l+1)=\frac{n}{2}\sum_{l=0}^{n}\psi(l+1)\psi(n-l+1). \quad \textbf{(B.55a)}$$

By reversing either or both series, (**B.56**) can be identified with several different combinations using the identities

$$\sum_{l=1}^{n}\psi(l+1)\psi(n-l+1)^2=\sum_{l=1}^{n}\psi(l+1)^2\psi(n-l+1)$$

and **(B.57)**

$$\sum_{l=1}^{n}\psi(l+1)\psi'(n-l+1)=\sum_{l=1}^{n}\psi'(l+1)\psi(n-l+1).$$

In a similar vein, consider the limit $\beta\to n, n>0$ in (**B.53**) and for the latter two cases employed previously, (the zero$^{th}$ order case reduces to (**B.16**)) obtain

$$\sum_{l=0}^{n-1}\frac{\psi(l+1)\psi(n-l)}{(n-l)}=\frac{1}{2}\psi(n+1)^3+\frac{\psi(n+1)}{2}(3\psi'(n+1)-\gamma^2-\frac{\pi^2}{2})+\frac{1}{2}\psi^{(2)}(n+1)+\zeta(3), \quad \textbf{(B.58)}$$

and

$$\sum_{l=0}^{n-1}\frac{\psi(l+1)\psi'(n-l)}{(n-l)}+\sum_{l=0}^{n-1}\frac{\psi(l+1)\psi(n-l)^2}{(n-l)}=\frac{1}{3}\psi(n+1)^4$$
$$+\psi(n+1)^2(2\psi'(n+1)-\frac{\pi^2}{6})+\psi'(n+1)^2-\frac{\pi^2}{6}\psi'(n+1) \tag{B.59}$$
$$+\frac{\psi(n+1)}{3}(4\psi^{(2)}(n+1)+8\zeta(3)+\frac{\gamma\pi^2}{2}+\gamma^3)+\frac{\psi^{(3)}(n+1)}{3}-\frac{\pi^4}{45}.$$

Again, by reversing each of the sums using the identities

$$\sum_{l=0}^{n-1}\frac{\psi(l+1)\psi(n-l)^2}{n-l}=\sum_{l=0}^{n-1}\frac{\psi(l+1)^2\psi(n-l)}{l+1}$$
and
$$\sum_{l=0}^{n-1}\frac{\psi(l+1)\psi'(n-l)}{n-l}=\sum_{l=0}^{n-1}\frac{\psi'(l+1)\psi(n-l)}{l+1} \tag{B.60}$$

a variety of sums can be identified. See also the note preceding (**B.47**).

Based on the hypergeometric identity **(16)** and its reversed form, that is

$$\sum_{l=0}^{n-1}\frac{\Gamma(l-a)}{\Gamma(l+1)(l-n)^2}=-\frac{\Gamma(n-a)}{\Gamma(n+1)}\Big(\sum_{l=0}^{n-1}\frac{\psi(2-n+a+l)}{(l+1)}-\psi(a+1-n)(\psi(n+1)+\gamma)\Big) \tag{B.61}$$

differentiating once with respect to $a$ gives

$$\sum_{l=0}^{n-1}\frac{\Gamma(l-a)\psi(l-a)}{\Gamma(l+1)(l-n)^2}=\frac{\Gamma(n-a)}{\Gamma(n+1)}\Big(-\psi(n-a)\sum_{l=0}^{n-1}\frac{\psi(2-n+a+l)}{(l+1)}+\sum_{l=0}^{n-1}\frac{\psi'(2-n+a+l)}{(l+1)}$$
$$+\big(\psi(a+1-n)\psi(n-a)-\psi'(a+1-n)\big)\big(\psi(n+1)+\gamma\big)\Big) \tag{B.62}$$

and in the limit $a\to -1$, we find

$$\sum_{l=0}^{n-1}\frac{\psi(l+1)}{(l-n)^2}=-\psi(n+1)\sum_{l=0}^{n-1}\frac{\psi(n-l)}{(l+1)}-\sum_{l=0}^{n-1}\frac{\psi'(n-l)}{(l+1)}+(\psi(1+n)^2+\psi'(1+n))(\psi(n+1)+\gamma). \tag{B.63}$$

In a sense, this generalizes (**B.3)**.

From the literature, Speiss [**30**, Lemma 18] (presented as a sum of harmonic numbers, rewritten here) provides the following transformation

$$\sum_{l=0}^{k}\frac{\psi(l+1)}{l+m+1}=-\sum_{l=1}^{m}\frac{\psi(l)}{l+k}+\frac{\gamma}{m+1}+\frac{\psi(k+1)}{l+m+1}+\tfrac{1}{2}(\psi'(k+1)+\psi'(m+1))-\frac{\pi^2}{6}$$
$$-\tfrac{1}{2}(\psi(k+1)+\psi(m+1))(\psi(k+1)+\psi(m+1)-2\psi(k+m+1)), \quad \textbf{(B.64)}$$

which, in conjunction with **(18)** using $b=1$ and $c=m+1$, produces the variation

$$\sum_{l=1}^{m}\frac{\psi(l)}{l+k}=\sum_{l=0}^{k}\frac{\psi(l+m+1)}{l+1}-\gamma\psi(m)+\psi(k+m+1)(\psi(m)-\tfrac{1}{k+1})-\frac{1}{2}(\psi(m)+\psi(k))^2$$
$$+\frac{1}{2}(\psi'(k)+\psi'(m))-\frac{\psi(m)+\psi(k+1)}{k}-\frac{\pi^2}{6}. \quad \textbf{(B.65)}$$

In the case that $k=0$ we obtain (also see **(B.13**) and **(B.20)**)

$$\sum_{l=1}^{m}\frac{\psi(l)}{l}=\frac{1}{2}(\psi(m)^2-\gamma^2+\psi'(m)-\frac{\pi^2}{6})+\frac{\psi(m)}{m}. \quad \textbf{(B.66)}$$

And, by considering its even and odd terms, we find the related alternating sum

$$\sum_{l=1}^{m}(-1)^l\frac{\psi(l)}{l}=-\tfrac{1}{2}\sum_{l=1}^{\lfloor\frac{m+1}{2}\rfloor}\frac{\psi(l)}{l-\frac{1}{2}}-\frac{(-1)^m}{2}[\psi(m)^2+\psi'(m)-\tfrac{1}{2}\psi'(\tfrac{m}{2})]+\psi(m)[\tfrac{(-1)^m+1}{2}\psi(\tfrac{m}{2})+\tfrac{1}{m}]$$
$$+\tfrac{(-1)^m-1}{2}[\psi(\tfrac{m}{2})(\ln(2)+\tfrac{1}{m})+\ln^2(2)+\tfrac{2}{m}\ln(2)]-\tfrac{1}{4}\psi(\tfrac{m}{2})^2+\frac{\pi^2}{24}-\frac{\gamma^2}{4}. \quad \textbf{(B.67)}$$

Finally, as discussed following **(69)** we find

$$\sum_{l=0}^{m-1}\frac{\psi(l+1)}{(l-m)^2}-\sum_{l=0}^{m-1}\frac{\psi(l+1)}{(l+1)^2}=(\tfrac{1}{6}\pi^2-\psi'(m+1))(\gamma+\psi(m+1))-\tfrac{1}{2}\psi^{(2)}(m+1)-\zeta(3). \quad \textbf{(B.68)}$$

See also **(B.6),** (**B.10**), **(B.48)** and **(B.63**).

## Appendix C: Infinite Sums

The following useful and related infinite sums, based on the methods of this paper, and needed (or obtained) elsewhere [**19**], are presented* with minimal proof.

$$\sum_{n_1=0}(\frac{1}{n_1+\frac{1}{2}})\sum_{n_2=0}(\frac{1}{n_{(2)}+\frac{3}{2}})(\frac{1}{n_{(2)}+1})^2=7\zeta(3)-\frac{2\pi^2}{3} \quad \textbf{(C.1)}$$

$$\sum_{n_1=0}(\frac{1}{n_1+\frac{1}{2}})\sum_{n_2=0}(\frac{1}{n_{(2)}+1})(\frac{1}{n_{(2)}+\frac{3}{2}})=2\sum_{l=0}\frac{\psi(l+\frac{3}{2})-\psi(l+1)}{l+\frac{1}{2}}=\frac{\pi^2}{3} \quad \textbf{(C.2)}$$

* $n_{(2)}=n_1+n_2$ ; $n_{(3)}=n_{(2)}+n_3$

$$\sum_{n_1=0}(\frac{1}{n_1+\frac{1}{2}})\sum_{n_2=0}(\frac{1}{n_{(2)}+1})(\frac{1}{n+n_{(2)}+\frac{5}{2}})=\frac{1}{(n+\frac{3}{2})}\Big(\sum_{l=0}^{n}\frac{\psi(\frac{3}{2}+l)-\psi(\frac{1}{2})}{(l+1)}+\frac{\pi^2}{6}\Big) \quad \textbf{(C.3)}$$

$$\sum_{l=0}\frac{\psi(l+\frac{1}{2})^2-\psi(l+1)^2+\psi'(l+\frac{1}{2})-\psi'(l+1)}{l+\frac{1}{2}}=-\frac{2\pi^2}{3}\psi(\tfrac{1}{2}) \quad \textbf{(C.4)}$$

$$\sum_{n_1=0}\frac{1}{(n_1+\frac{1}{2})}\sum_{n_2=0}\frac{1}{(n_1+n_2+1)}[\psi(3+n_{(2)})-\psi(\tfrac{5}{2}+n_{(2)})]$$
$$=\frac{1}{2}\sum_{n_1=0}(\frac{1}{n_1+\frac{1}{2}})\sum_{n_2=0}(\frac{1}{n_{(2)}+1})\sum_{n_3=0}(\frac{1}{n_{(3)}+\frac{5}{2}})(\frac{1}{n_{(3)}+3})=4\log(2)+\tfrac{3}{2}\zeta(3)-\frac{\pi^2}{3} \quad \textbf{(C.5)}$$

Evaluate $p \to q$ in **(37)**:

$$\sum_{l=0}\frac{\psi(l+m+1)}{(l+q)^2}=\psi(q)\psi'(q-m)+\psi'(q)\psi(q-m)$$
$$-\tfrac{1}{2}\psi^{(2)}(q)-\psi(q)\psi'(q)-\sum_{l=0}^{m-1}\frac{\psi(l+1)}{(q-l-1)^2} \quad \textbf{(C.6)}$$

then if $q\to m$:

$$\sum_{l=0}\frac{\psi(l+m+1)}{(l+m)^2}=\psi(m)(\frac{\pi^2}{6}-\psi'(m))-\psi^{(2)}(m)-\gamma\psi'(m)-\sum_{l=0}^{m-2}\frac{\psi(l+1)}{(l-m+1)^2} \quad \textbf{(C.7)}$$

$$\sum_{l=0}[\frac{\psi(a+l+n)}{(a+l)}-\frac{\psi(a+l+q+n)}{(a+l+q)}]$$
$$=\Big\{\frac{1}{2}(\psi'(a)-\psi(a)^2)-\psi(a)[\psi(n)-\psi(1)]+\sum_{l=0}^{n-2}\frac{\psi(a+l+1)}{(l+1)}\Big\}-\{a\to a+q\} \quad \textbf{(C.8)}$$

Generalizations of **(56)**:

$$\sum_{l=0}[\frac{\psi(l+k+j+\frac{1}{2})-\psi(l+k+1)}{l+k+\frac{1}{2}}]$$
$$=\frac{1}{2}\psi(\tfrac{1}{2})^2-\frac{\pi^2}{12}+\psi'(k+\tfrac{1}{2})(\tfrac{1}{2}-\delta_{j0})-\tfrac{1}{2}\psi(k+\tfrac{1}{2})^2+\sum_{l=0}^{k-1}\frac{\psi(l+1)}{(l+\frac{1}{2})}$$
$$+\sum_{l=0}^{j-2}\frac{\psi(k+l+\frac{3}{2})-\psi(k+\frac{1}{2})}{(l+1)} \qquad k>0 \quad \textbf{(C.9)}$$

$$\sum_{l=0}\frac{[\psi(l+N+1)-\psi(l+j+1)]}{(l+y)}$$
$$=\sum_{k=0}^{N-j-1}\frac{\psi(j+1+k)}{(j+1+k-y)}-\psi(y)\big(\psi(N+1-y)-\psi(j+1-y)\big),\ N>j \quad \textbf{(C.10)}$$

$$\sum_{l=0}\frac{[\psi(l+x+N+1)-\psi(l+x+1)]}{(l+y)}=\sum_{k=0}^{N-1}\frac{\psi(x+k+1)-\psi(y)}{(x+k+1-y)} \quad \textbf{(C.11)}$$

$$x, y \neq 0, -1, \cdots$$

$$\sum_{l=0} \frac{[\psi(l+x+N+1)-\psi(l+x+1)]}{(l+x)} = \sum_{k=0}^{N-1} \frac{\psi(x+k+1)-\psi(x)}{(k+1)} \quad \textbf{(C.12)}$$

$$x \neq -1, -2, \cdots$$

From **(62)** and [**13**, Eq. 9)[*]]

$$\sum_{l=0} \frac{\psi'(\frac{1}{2}+l)}{(\frac{1}{2}+l)} = \tfrac{7}{2}\zeta(3) + \pi^2 \log(2). \quad \textbf{(C.13)}$$

Evaluate the limit $k \to \infty$ in **(B.9)** to obtain a variation of **(66)**

$$\sum_{l=0} \frac{\psi^{(n)}(b+l)}{(b+l)^{n+1}} = \frac{(-1)^{n+1}}{2}\left\{\frac{1}{n!}\psi^{(n)}(b)^2 + \frac{n!}{(2n+1)!}\psi^{(2n+1)}(b)\right\},\ n > 0. \quad \textbf{(C.14)}$$

See the note following **(60)** for the (divergent) case n=0.

In (**B.12**) when $k \to \infty$ we find

$$\sum_{l=0} \frac{\psi'(1+l)}{(1+l)^2(2+l)^2} = \frac{7\pi^4}{180} + \frac{\pi^2}{6} - 5. \quad \textbf{(C.15)}$$

From **[19, Eqs. (5.1.7), (5.1.10) and (2.10)]]** we have

$$\sum_{l=0} \frac{t_l(1)}{(l+1)(l+\frac{3}{2})(l+2)} = \tfrac{2}{3}\pi^2 - 8\log(2) \quad \textbf{(C.16)}$$

and

$$\sum_{l=0} \frac{t_l(2)}{(l+\frac{3}{2})(l+\frac{5}{2})(l+3)} = \frac{2}{3}\pi^2 - \frac{16}{3}\log(2) - 2\zeta(3) \quad \textbf{(C.17)}$$

where [**20**]

$$t_l(1) = \psi(l+\tfrac{3}{2}) - \psi(\tfrac{1}{2}) = H_{l+1}(\tfrac{1}{2})$$
$$t_l(2) = \sum_{k=0}^{l} \frac{\psi(k+\frac{3}{2}) - \psi(\frac{1}{2})}{(k+1)}\ . \quad \textbf{(C.18)}$$

From the Thomae transformations [**18**, section 3.13.3] as employed in **(51)**, variations of **(52)** are:

$$\sum_{l=0} \left( \frac{\psi'(l+1)}{(q+l)} - \frac{\psi(q+1+l)}{(l+1)^2} \right) = -\frac{\pi^2}{6}\psi(q) \quad \textbf{(C.19)}$$

---

[*] which reads: $\sum_{l=0} \frac{\psi(l+\frac{1}{2})}{(l+\frac{1}{2})^2} = -(7\zeta(3)+\gamma\pi^2)/2$ - see (**C.28**).

$$\sum_{l=0}\frac{\psi'(l+1)}{(q+l)}=\gamma\psi'(q)+\psi(q)\psi'(q)-\tfrac{1}{2}\psi^{(2)}(q)+{}_4F_3\left({1,1,1,1 \atop 2,2,q+1}\middle|1\right)/q$$
$$=2\zeta(3)-\frac{\pi^2}{6}(\gamma+\psi(q))-\frac{1}{\Gamma(1-q)}\sum_{l=0}\frac{\Gamma(2-q+l)\psi'(l+2)}{\Gamma(2+l)(1+l)} \quad \textbf{(C.20)}$$

$$\sum_{l=0}\left(\frac{\psi'(l+1)}{(q+l)}-\frac{\psi(q+l)}{(q+l)^2}\right)=\gamma\psi'(q)-\tfrac{1}{2}\psi^{(2)}(q) \quad \textbf{(C.21)}$$

$$\sum_{l=0}\left(\frac{\psi'(l+1)-\psi'(l+q)}{(q+l)}\right)=\gamma\psi'(q)+\psi(q)\psi'(q) \quad \textbf{(C.22)}$$

$$\sum_{l=0}\left(\frac{\psi'(l+1)}{(q+l)}+\frac{\psi(2-q+l)}{(l+1)^2}\right)=\frac{\pi^2}{6}\psi(q)+\pi\cot(\pi q)\left(\frac{\pi^2}{6}-\psi'(q)\right)-\tfrac{1}{2}\psi^{(2)}(q) \quad \textbf{(C.23)}$$

Therefore, comparing **(C.19)** and (**C.23**) at $q=\frac{1}{2}$ gives

$$\sum_{l=0}\frac{\psi(\frac{1}{2}+l+n)}{(1+l)^2}=\left(\frac{\pi^2}{6}+\frac{1}{(n-\frac{1}{2})^2}\right)\psi(n-\tfrac{1}{2})+\gamma\psi'(n-\tfrac{1}{2})-\tfrac{1}{2}\psi^{(2)}(n-\tfrac{1}{2})$$
$$-\tfrac{1}{2}\pi^2\gamma-\tfrac{7}{2}\zeta(3)-\sum_{l=0}^{n-1}\frac{\psi(\frac{1}{2}+l)}{(\frac{1}{2}+l)^2}, \quad \textbf{(C.24)}$$

$$\sum_{l=0}\frac{\psi'(1+l)}{(\frac{1}{2}+l)}=\tfrac{7}{2}\zeta(3), \quad \textbf{(C.25)}$$

and

$$\sum_{l=0}\frac{\psi(\frac{1}{2}+l)}{(1+l)^2}=\pi^2\left(\frac{1}{3}-\frac{\gamma}{6}-\frac{\ln(2)}{3}\right)-8\ln(2)+\tfrac{7}{2}\zeta(3)\cdot \quad \textbf{(C.26)}$$

Equation **(C.24**) generalizes **[9, Eq.(15)]**. Another more general form of the above appears in **[13, Eq. (8a)]**, that is

$$\sum_{l=1}\frac{\psi(\frac{1}{2}+l)}{l^{2n}}=\psi(\tfrac{1}{2})\zeta(2n)+\tfrac{1}{2}\zeta(2n+1)(2^{2n+1}-1)$$
$$-\sum_{l=1}^{n-1}(2^{2n+1-2l}-1)\zeta(2l)\zeta(2n+1-2l), \quad \textbf{(C.27)}$$

along with its related form

$$\sum_{l=1}\frac{\psi(l-\frac{1}{2})}{(l-\frac{1}{2})^{2n}}=-\gamma\zeta(2n)(2^{2n}-1)-\tfrac{1}{2}\zeta(2n+1)(2^{2n+1}-1)$$
$$-\sum_{l=1}^{n-1}(2^{2l}-1)\zeta(2l)\zeta(2n+1-2l). \quad \textbf{(C.28)}$$

Applying the limit $k\to\infty$ to **(21)** gives

$$\sum_{l=0}\left(\frac{\psi'(l+b)}{(c+l)}-\frac{\psi(c+l)}{(l+b)^2}\right)=-\psi'(b)\psi(c)+(\psi(b)-\psi(c))/(b-c)^2-\psi'(b)/(b-c). \qquad \textbf{(C.29)}$$

For the limit $b\to c$ see **(62)**. After operating on **(C.29)** with the operator $\frac{d}{dc}$ we obtain

$$\sum_{l=0}\left(\frac{\psi'(l+b)}{(c+l)^2}+\frac{\psi'(c+l)}{(l+b)^2}\right)=\psi'(b)(\psi'(c)+1/(b-c)^2)+2(\psi(c)-\psi(b))/(b-c)^3 +\psi'(c)/(b-c)^2 \qquad \textbf{(C.30)}$$

and, in the limit $c\to b$

$$\sum_{l=0}\frac{\psi'(l+b)}{(b+l)^2}=\tfrac{1}{2}\psi'(b)^2+\tfrac{1}{12}\psi^{(3)}(b) \qquad \textbf{(C.31)}$$

From [**19**] and resolving* [**12**, Eq. (55.6.4)] find:

$$\sum_{l=0}\frac{x^{2l+1}[\psi(1+l)-\psi(1)]}{(l+\frac{1}{2})}=\frac{\pi 2}{6}-\log^2(\tfrac{1+x}{2})-2Li_2(\tfrac{1-x}{2})+\tfrac{1}{2}\log^2(\tfrac{1+x}{1-x})-\log(\tfrac{1+x}{1-x})\log 4 \qquad \textbf{(C.32)}$$

Operate on **(C.32)** with $\frac{d}{dx}$ to obtain

$$\sum_{l=0}x^l(\psi(1+l)-\psi(1))=\frac{\ln(1-x)}{(x-1)} \qquad \textbf{(C.33)}$$

reproducing a known result [**12**, Eq. (55.3.1)] ( $x\neq 1$ ).

From [**21**, Eqs. (9), (10) and (12)] we have:

$$\sum_{l=0}\frac{\psi(q+\frac{1}{2}l)-\psi(q+\frac{1}{2}l+\frac{1}{2})}{(1+l-2q)}=\frac{(\psi(q+\frac{1}{2})-\psi(q))^2}{4}+(\psi(2q)-\psi(2q-\tfrac{1}{2}))(\psi(1-2q)-\psi(2q)) \qquad \textbf{(C.34)}$$

equivalent to

$$\sum_{l=0}\frac{\psi(q+l)-\psi(q+l+\frac{1}{2})}{(1+l-q)(\frac{1}{2}+l-q)}=(\psi(q+\tfrac{1}{2})-\psi(q))^2+4\pi\cot(2\pi q)(\psi(2q)-\psi(2q-\tfrac{1}{2}))-\frac{\pi\cot(\pi q)}{(q-\frac{1}{2})} \qquad \textbf{(C.35)}$$

and

$$\sum_{l=0}\frac{(-1)^l\psi(q+l)}{(l+2q-1)}=\frac{\psi(q)(\psi(q)-\psi(q-\frac{1}{2}))}{2}-\frac{(\psi(\frac{1+q}{2})-\psi(\frac{q}{2}))^2}{8}. \qquad \textbf{(C.36)}$$

The result **(C.34**) generalizes [**9**, Eq(24)], which corresponds to the case q=0. If $q=\frac{1}{2}$(equivalent to [**12**, **Eq. (55.6.2)** when $x=i$]) then

$$\sum_{l=0}\frac{(-1)^l\psi(\frac{1}{2}+l)}{(l+1)}=\ln(2)(-\gamma-2\ln(2))+\frac{\pi^2}{8}-\pi+2\ln(2), \qquad \textbf{(C.37)}$$

which, by grouping terms in pairs can equivalently be written as

$$\tfrac{1}{4}\sum_{l=0}\frac{\psi(\tfrac{1}{2}+2l)}{(l+1)(\tfrac{1}{2}+l)}=\ln(2)(-\gamma-2\ln(2))+\frac{\pi^2}{8}-\frac{5\pi}{6}+3\ln(2). \qquad \textbf{(C.38)}$$

---

* notice that [**12**,Eq. (55.6.7)] is incorrect – omit the term "*2x*";

Applying the limit $k\to\infty$ to (**28**) with $a-b\geq 2$ gives

$$\sum_{l=0}\frac{\Gamma(b+l)}{\Gamma(a+l)}\psi(b+l)=\frac{\Gamma(b)\psi(b)}{(a-b-1)\Gamma(a-1)}+\frac{\Gamma(b)}{\Gamma(a-1)(a-b-1)^2}, \qquad \textbf{(C.39)}$$

yielding the convergent series

$$\sum_{l=0}\frac{\psi(b+l)}{(b+l)(b+l+1)}=\frac{1}{b}(\psi(b)+1) \qquad \textbf{(C.40)}$$

when $a=b+2$, generalizing [**17**, corollary 2.3]. To **(C.39)** apply the ordered operators $\lim_{a\to b+2}\frac{\partial}{\partial b}$ and $\lim_{a\to b+2}\frac{\partial}{\partial a}$ independently, then add the two results to eventually obtain

$$\sum_{l=0}\frac{\psi'(b+l)}{(b+l)(b+l+1)}=\frac{1}{b}(\psi'(b)+1)-\psi'(b) \qquad \textbf{(C.41)}$$

and

$$\sum_{l=0}\psi(b+l)\Big(\frac{1}{(b+1+l)^2}-\frac{1}{(b+l)^2}\Big)=-\frac{1}{b}(\psi(b)/b+1)+\psi'(b+1). \qquad \textbf{(C.42)}$$

Useful sums deduced from among all the previously compiled entries follow:

$$2\sum_{l=1}\Big(\frac{\psi(l+(j+1)/2)}{l+(j-1)/2}-\frac{\psi(l+j/2)}{l+j/2}\Big)=\psi(j/2)^2+\psi'(j/2)-\psi((j+1)/2)^2+\psi'((j+1)/2)+4\psi(j/2)/j \qquad \textbf{(C.43)}$$

$$\sum_{l=1}\Big(\frac{\psi(l+j-\frac{1}{2})}{l+j-1}-\frac{\psi(l+j)}{l+j-\frac{1}{2}}\Big)=\frac{\psi(j-1)}{j-\frac{1}{2}}-\psi(j)\psi(j-\tfrac{1}{2})+\frac{2}{j-1}-\frac{2}{j-\frac{1}{2}}-\frac{\pi^2}{6}+(\gamma+2\ln(2))^2+2\sum_{l=0}^{j-2}\frac{\psi(l+1)}{l+\frac{1}{2}} \qquad \textbf{(C.44)}$$

$$\sum_{l=1}\Big(\frac{\psi(l+j)}{l}-\frac{\psi(l+j-\frac{1}{2})}{l-\frac{1}{2}}\Big)=\frac{1}{2}(\psi(j)^2-\psi'(j))-2\psi(j-\tfrac{1}{2})+(2-2\ln(2))\psi(j)-\frac{\pi^2}{12}-2-(\gamma-2+2\ln(2))(\gamma+2-2\ln(2))/2-\sum_{l=0}^{j-2}\frac{\psi(l+\frac{1}{2})}{l+1} \qquad \textbf{(C.45)}$$

and, by solving a single step recursion equation

$$\sum_{l=1}\Big(\frac{\psi(l+j)}{l-\frac{1}{2}}-\frac{\psi(l+j-\frac{1}{2})}{l}\Big)=2\psi(j)-\frac{1}{2}\big(\psi(\tfrac{1}{2}+j)^2-\psi'(j+\tfrac{1}{2})\big)+2\ln^2(2)$$
$$+(2\ln(2)+\frac{1}{j-\frac{1}{2}}-2)\psi(\tfrac{1}{2}+j)+\gamma^2/2+(2\gamma-4)\ln(2) \tag{C.46}$$
$$-(2-\frac{1}{j-\frac{1}{2}})\gamma+4-\frac{\pi^2}{12}+\sum_{l=0}^{j-2}\frac{\psi(l+1)}{l+\frac{3}{2}}\,.$$

By considering the odd and even terms of the sum $\sum_{l=1}\frac{(-1)^l\psi(l+n)}{l(l+n)}$ individually and applying the above, as well as simple results summarized in Appendix E, the following two results for odd and even values of $n$ respectively can be found:

$$(2j-1)\sum_{l=1}\frac{(-1)^l\psi(l+2j-1)}{l(l+2j-1)}=2\ln^2(2)+(\psi(j-1/2)-2\psi(2j-1)-2)\ln(2)-\psi(j-1/2)$$
$$+\psi(j)-(\gamma^2+\psi'(j)-\psi(j-1/2)^2)/4-\frac{\pi^2}{24}+\frac{\psi(2j-1)}{(2j-1)}-\frac{1}{2}\sum_{l=0}^{j-2}\frac{\psi(l+1/2)}{l+1} \tag{C.47}$$

and

$$j\sum_{l=1}\frac{(-1)^l\psi(l+2j)}{l(l+2j)}=-\ln^2(2)+(1-\psi(j+\tfrac{1}{2})/2-\gamma/2)\ln(2)+\psi(j+\tfrac{1}{2})/2$$
$$-\big(\psi'(j+\tfrac{1}{2})-\psi(j)^2+\gamma^2-\pi^2/6\big)/8-\psi(j)/2-1+\gamma/2 \tag{C.48}$$
$$+\big(\psi(j+\tfrac{1}{2})+\psi(j)+2\ln(2)\big)/(8j)-\frac{1}{4}\sum_{l=0}^{j-2}\frac{\psi(l+1)}{l+\frac{3}{2}}\,.$$

From an unpublished source **[11]**, the following was contributed

$$\sum_{l=1}(-1)^l\frac{\psi(1+\frac{1}{2}l)-\psi(\frac{1}{2}+\frac{1}{2}l)}{l^2}=\tfrac{1}{2}\zeta(3)-\frac{\pi^2\ln(2)}{6} \tag{C.49}$$

and

$$\sum_{m=1}\sum_{l=1}(-1)^{l+m}\frac{\psi(l+m)}{lm(l+m)}=\frac{11}{240}\pi^4-\frac{\gamma}{4}\zeta(3)-4Li_4(\tfrac{1}{2})-\frac{7}{2}\zeta(3)\ln(2)-\frac{1}{6}(\ln^4(2)-\pi^2\ln^2(2))\,, \tag{C.50}$$

the first of which (see also **(C.36)** and **(C.37)**) can be proven using the results contained previously herein by considering its even and odd terms. The inner sum(s) of the second is given by **(C.47)** and **(C.48);** although **(C.50)** was offered without proof, it is numerically verifiable, and, temporarily accepting its validity, a number of other interesting and useful results can be deduced. Adding the even and odd terms of the inner sum of **(C.50)** and utilizing many of the results found in Appendix E, eventually gives

$$\sum_{n=1}\frac{1}{(n+\frac{1}{2})^2}\sum_{m=1}^{n}\frac{\psi(1+m)}{(\frac{1}{2}+m)}-\sum_{n=1}\frac{1}{(n+1)^2}\sum_{m=1}^{n}\frac{\psi(\frac{1}{2}+m)}{(1+m)}=(\tfrac{3}{2}\psi(\tfrac{1}{2})-7\ln(2)+3)\zeta(3)$$
$$+((\ln(2)-\tfrac{5}{6})\psi(\tfrac{1}{2})-2\ln(2)+\tfrac{2}{3}\ln^2(2))\pi^2-\tfrac{2}{3}\ln^4(2)+\frac{83}{360}\pi^4-16Li_4(\tfrac{1}{2}) \tag{C.51}$$

or, after employing **(23)**

$$\sum_{n=1}\left(\frac{1}{(n+\frac{1}{2})^2}+\frac{1}{(n+1)^2}\right)\sum_{m=1}^{n}\frac{\psi(m+1)}{(\frac{1}{2}+m)}=(7\psi(\tfrac{1}{2})+7\ln(2))\zeta(3)$$
$$+((\ln(2)-\tfrac{4}{3})\psi(\tfrac{1}{2})-\tfrac{8}{3}\ln(2)+\tfrac{1}{3}\ln^2(2))\pi^2-\tfrac{1}{3}\ln^4(2)+\frac{79}{360}\pi^4-8Li_4(\tfrac{1}{2}). \qquad \textbf{(C.52)}$$

Subtract **(C.51)** from **(C.52)** to obtain

$$\sum_{n=1}\frac{1}{(n+1)^2}\sum_{m=1}^{n}\left(\frac{\psi(m+\frac{1}{2})}{(1+m)}+\frac{\psi(m+1)}{(\frac{1}{2}+m)}\right)=(\tfrac{11}{2}\psi(\tfrac{1}{2})+14\ln(2)-3)\zeta(3)$$
$$-(\tfrac{1}{2}\psi(\tfrac{1}{2})+\tfrac{2}{3}\ln(2)+\tfrac{1}{3}\ln^2(2))\pi^2+\tfrac{1}{3}\ln^4(2)-\tfrac{1}{90}\pi^4+8Li_4(\tfrac{1}{2}), \qquad \textbf{(C.53)}$$

Again apply **(23)** to the inner sum, and notice that all the sums so-created appear in Appendix E. After considerable computation and simplification, it is found that **(C.53)** reduces to an identity, thereby verifying **(C.50).**

Alternatively, reorder **(C.52)** along a diagonal of the summation grid, interchange the summations, simplify using results from Appendix E and eventually obtain

$$\sum_{l=1}\frac{\psi(l)\psi'(2l-1)}{(\frac{1}{2}+l)}=\tfrac{7}{4}(\psi(\tfrac{1}{2})+\ln(2)+1)\zeta(3)-\tfrac{1}{12}\ln^4(2)+\frac{79}{1440}\pi^4+2\ln^2(2)-5$$
$$+(\tfrac{1}{4}(\ln(2)-\tfrac{1}{2})\pi^2+2\ln(2)-\tfrac{5}{2})\psi(\tfrac{1}{2})+(-\tfrac{1}{2}+\tfrac{1}{2}\ln(2)+\tfrac{1}{12}\ln^2(2))\pi^2-2Li_4(\tfrac{1}{2}). \qquad \textbf{(C.54)}$$

## Appendix D (A Special Case)

The following is the computer output defining the result associated with the generalization of **(B.43)** after differentiating with respect to the variable β. In the following, the variable $B\equiv\beta-m$ and the various parameters $Y_n$ are defined following the main result. For the special cases $\beta=m$ and $\beta=m+1$, see (**B.45**) and (**B.46**).

$$Y_1 B\left(\sum_{k=0}^{m}\left(-\frac{(-1)^k\Psi(k+1)^2\Psi(k-\beta+1)}{\Gamma(m-k+1)\Gamma(k-\beta+1)}\right)\right)=\left(Y_2-\frac{1}{B}\right)\left(\sum_{k=0}^{m}\frac{\Psi(m-1-k-\beta)}{k+1}\right)+\sum_{k=0}^{m}\left(-\frac{\Psi^{(1)}(m-1-k-\beta)}{k+1}\right)+\left(Y_2-\frac{1}{B}\right)\gamma^2+\left(Y_3-\frac{Y_4}{B}\right)\gamma+\frac{Y_5}{2}-\frac{Y_6}{2B}$$

$$Y_1=\Gamma(m+1)\Gamma(-\beta)$$

$$Y_2=\Psi(-\beta)$$

$$Y_3=-2\Psi(-\beta)^2+\Psi(-\beta)\Psi(-B)+2\Psi^{(1)}(-\beta)-\Psi^{(1)}(-B)$$

$$Y_4=-2\Psi(-\beta)+\Psi(-B)$$

$$Y_5=\Psi(-\beta)\Psi(-B)^2-2\left(\Psi(-\beta)^2+\Psi(-\beta)\Psi(m+1)-\Psi^{(1)}(-\beta)+\Psi^{(1)}(-B)\right)\Psi(-B)+\Psi(-\beta)^3-\frac{2\Psi(-\beta)^2}{m+1}+\left(3\Psi^{(1)}(-B)-3\Psi^{(1)}(-\beta)\right)\Psi(-\beta)+2\Psi^{(1)}(-B)\Psi(m+1)+\frac{2\Psi^{(1)}(-\beta)}{m+1}-\Psi^{(2)}(-B)+\Psi^{(2)}(-\beta)+\frac{2}{(1+\beta)^2(m+1)}$$

$$Y_6=-2\Psi(-B)\Psi(m+1)+\Psi(-B)^2-2\Psi(-\beta)\Psi(-B)+\Psi(-\beta)^2-\frac{2\Psi(-\beta)}{1+\beta}+\Psi^{(1)}(-B)-\Psi^{(1)}(-\beta)-\frac{2}{(m+1)(1+\beta)}$$

## Appendix E (Simple Special Cases)

The following are commonly needed, useful and interesting cases that can be derived from the more general entries listed above, or extracted from the literature.

$$\sum_{l=0}\frac{\psi(l+1/2)}{(l+1/2)^2}=-\frac{7}{2}\zeta(3)-\pi^2\gamma/2 \qquad \textbf{(E.1)}^{*}$$

$$\sum_{l=0}\frac{\psi(l+1)}{(l+1/2)^2}=-\frac{\pi^2}{2}(\gamma+2\ln(2))+7\zeta(3) \qquad \textbf{(E.2)}$$

$$\sum_{l=0}\frac{\psi(l+1)}{(l+1)^2}=-\frac{1}{6}\pi^2\gamma+\zeta(3) \qquad \textbf{(E.3)}$$

$$\sum_{l=0}\frac{\psi(l+1)}{(l+1)^3}=\frac{\pi^4}{360}-\zeta(3)\gamma \qquad \textbf{(E.4)}^{\dagger}$$

$$\sum_{l=0}\frac{\psi(l+1)}{(l+\frac{1}{2})^3}=-7(2\ln(2)+\gamma)\zeta(3)+\frac{\pi^4}{8} \qquad \textbf{(E.5)}$$

$$\sum_{l=0}\frac{\psi(l+\frac{3}{2})}{(l+1)^2}=-\tfrac{1}{3}\pi^2(\gamma/2+\ln(2))+\tfrac{7}{2}\zeta(3) \qquad \textbf{(E.6)}^{\ddagger}$$

$$\sum_{l=0}\frac{\psi'(l+1)}{(l+\frac{1}{2})}=\frac{7}{2}\zeta(3) \qquad \textbf{(E.7)}^{\S}$$

$$\sum_{l=0}\frac{\psi'(l+1)}{(l+\frac{1}{2})^2}=-28\ln(2)\zeta(3)+\tfrac{151}{360}\pi^4-\tfrac{4}{3}\ln^4(2)+\tfrac{4}{3}\pi^2\ln^2(2)-32Li_4(\tfrac{1}{2}) \qquad \textbf{(E.8)}$$

$$\sum_{l=0}\frac{\psi'(l+\frac{1}{2})}{(l+\frac{1}{2})^2}=\frac{5\pi^4}{24} \qquad \textbf{(E.9)}^{**}$$

$$\sum_{l=0}\frac{\psi(l+1)^2}{(l+1)^2}=\frac{11\pi^4}{360}+\frac{\pi^2\gamma^2}{6}-2\zeta(3)\gamma \qquad \textbf{(E.10)}^{\dagger\dagger}$$

$$\sum_{l=0}\frac{\psi'(l+\frac{1}{2})}{(l+\frac{1}{2})}=\pi^2\ln(2)+\frac{7}{2}\zeta(3) \qquad \textbf{(E.11)}^{\ddagger\ddagger}$$

---

[*] Repeating footnote to **(C.13).**

[†] See **(72).**

[‡] Repeating **(C.32).** See also **[9, Eq(15)]**.

[§] Repeating **(C.25).**

[**] Special case of (**C.31**)

[††] See **(72)** and **[9, Eq.(9)].**

[‡‡] Repeating **(C.13).**

$$\sum_{l=1}\left(\frac{\psi(l+1)}{l}-\frac{\psi(l+\frac{1}{2})}{l+\frac{1}{2}}\right)=\tfrac{1}{2}(\gamma+2\ln(2))^2+\pi^2/3-\gamma^2/2-2\gamma-4\ln(2)$$ **(E.12)*** 

$$\sum_{l=1}\left(\frac{\psi(l+\frac{1}{2})}{l}-\frac{\psi(l+1)}{l+\frac{1}{2}}\right)=-2\gamma(1-\ln(2))-\pi^2/6+4\ln^2(2)$$ **(E.13)†**

$$\sum_{l=0}\frac{\psi'(l+\frac{1}{2})}{(l+1)^2}+\frac{\psi'(l+1)}{(l+\frac{1}{2})^2}=-32\ln(2)+\frac{8\pi^2}{3}+\frac{\pi^4}{12}$$ **(E.14)‡**

$$\sum_{l=0}\frac{\psi(l+1)}{(l+1)(l+\frac{1}{2})}=\tfrac{1}{3}\pi^2-4\gamma\ln(2)-4\ln^2(2)$$ **(E.15)**

$$\sum_{l=1}\frac{\psi(l+\frac{1}{2})}{l(l+\frac{1}{2})}=\tfrac{1}{3}\pi^2+4(\gamma-2)\ln(2)+8\ln^2(2)-4\gamma$$ **(E.16)§**

$$\sum_{l=0}(-1)^l\frac{\psi(l+1)}{(l+1)}=-\tfrac{1}{2}\ln^2(2)-\gamma\ln(2).$$ **(E.17)**

$$\sum_{l=1}\frac{(-1)^l\psi(\frac{1}{2}+l)}{l}=\ln(2)(\gamma+2\ln(2))-\frac{\pi^2}{8}.$$ **(E.18)****

$$\sum_{l=1}(-1)^l\frac{\psi(\frac{1}{2}l)}{l^2}=\tfrac{1}{12}\gamma\pi^2+\frac{9}{8}\zeta(3)$$ **(E.19)**

$$\sum_{l=1}(-1)^l\frac{\psi'(\frac{1}{2}l)}{l^2}=-\tfrac{17}{360}\pi^4$$ **(E.20)**

The following are useful results extracted from **[2]**, **[7], [9], [10]**, [**25**] and/or [**37**]. From the general result [**25**, Eq. (3.11)] we have

$$\sum_{l=1}(-1)^l\frac{\psi(l+1)}{l^2}=\tfrac{1}{12}\gamma\pi^2-\tfrac{5}{8}\zeta(3).$$ **(E.21)**

Quoting **[9, Eq. (21)]**

$$\sum_{l=1}\frac{\psi(\frac{1}{2}\pm l)^2}{l}=-\tfrac{7}{4}\zeta(3)+\ln(2)(\tfrac{1}{2}\pi^2-\gamma^2)+(\tfrac{1}{4}\pi^2-4\ln^2(2))\gamma-4\ln^3(2),$$ **(E.22)**

and from [**25**, **Eq. (4.12)]** or **[7, Eq. (38)]**

$$\sum_{l=0}\frac{\psi'(l+1)}{l+1}=2\zeta(3),$$ **(E.23a)**

$$\sum_{l=0}(-1)^l\frac{\psi'(l+1)}{l+1}=-\tfrac{1}{4}\zeta(3)+\tfrac{1}{4}\pi^2\ln(2),$$ **(E.23b)**

along with, as a consequence of **[10**, **Theorem 3.1]** or **(C.14),**

---

* Special case of **(C.43)**
† Special case of **(C.44)**.
‡ See (**E.8**) **E.8**
§ From [**37**], second unnumbered equation following **(1.9)**, corrected for an incorrect exponent minus sign
**See **(C.37)**

$$\sum_{l=1}\frac{\psi'(l+1)}{l^2}=\tfrac{1}{120}\pi^4. \quad \textbf{(E.24)}$$

From corrected **[9, Eq. (22)],** (apply a factor of two to the right-hand side; also see **(E.6)**)

$$\sum_{l=1}\frac{\psi(\frac{1}{2}\pm l)^2}{l^2}=\pi^4/8+\psi(\tfrac{1}{2})^2\pi^2/6+7\psi(\tfrac{1}{2})\zeta(3). \quad \textbf{(E.25)}$$

Similarly, from the sixth (corrected) Example (3.6) of **[37]** – (remove the coefficient $\frac{7}{2}$; it should be unity) we obtain

$$\sum_{l=1}\frac{\psi(1+l)^2}{(l+\frac{1}{2})^2}=-\frac{61}{360}\pi^4+\pi^2(\tfrac{2}{3}\ln^2(2)+2\gamma\ln(2)+\tfrac{1}{2}\gamma^2)+\tfrac{4}{3}\ln^4(2)-14\gamma\zeta(3)-4\gamma^2+32Li_4(\tfrac{1}{2}). \quad \textbf{(E.26)}$$

From corrected **[9, Eq. (19)]**, (apply a minus sign to the right-hand side) or [**25**, **Eq. (3.11)** with *k=2*]

$$\sum_{l=1}(-1)^l\frac{\psi(\frac{1}{2}\pm l)}{l^2}=7\zeta(3)/2-2\pi G-\psi(\tfrac{1}{2})\pi^2/12. \quad \textbf{(E.27)}$$

From corrected **[9, Eq. (12)]**, (replace ½ ζ(3) with ¼ ζ(3); also in **[9, Eq. (11)]** replace $\frac{\pi^2}{12}$ with $\frac{\pi^2}{12}\ln(2)$), and **(E.17)**

$$\sum_{l=0}(-1)^l\frac{\psi(1+l)^2}{l+1}=\tfrac{1}{3}\ln^3(2)+\gamma\ln^2(2)-(\pi^2/12-\gamma^2)\ln(2)+\tfrac{1}{4}\zeta(3). \quad \textbf{(E.28)}$$

From corrected **[9, Eq. (8)],** (to both Eq. (8) and previous unnumbered equation, replace $\frac{1}{12}\pi^2\ln(2)$ with $\frac{1}{12}\pi^2\ln^2(2)$) or **[7, Eq.(30)]** we have

$$\sum_{l=1}(-1)^l\frac{\psi(l)}{l^3}=-\frac{\pi^4}{48}+2Li_4(\tfrac{1}{2})+\tfrac{7}{4}\ln(2)\zeta(3)-\tfrac{1}{12}\ln^2(2)\pi^2+\tfrac{1}{12}\ln^4(2)+\tfrac{3}{4}\gamma\zeta(3), \quad \textbf{(E.29)}$$

and from corrected **[9, Eq. (13)],** (contains the same error as **[9, Eq. (8)]** – see above)

$$\sum_{l=1}(-1)^l\frac{\psi(l)^2}{l^2}=\frac{11}{480}\pi^4-2Li_4(\tfrac{1}{2})-\tfrac{1}{12}\ln^4(2)-\tfrac{1}{4}(7\ln(2)+\gamma)\zeta(3)+\tfrac{\pi^2}{12}(\ln^2(2)-\gamma^2). \quad \textbf{(E.30)}$$

By considering the even and odd terms of **(E.29)** we obtain

$$\sum_{l=1}\frac{\psi(l+\frac{1}{2})}{l^3}=\zeta(3)(14\ln(2)+\psi(\tfrac{1}{2}))-\tfrac{2}{3}\pi^2\ln^2(2)+\tfrac{2}{3}\ln^4(2)-\tfrac{53}{360}\pi^4+16Li_4(\tfrac{1}{2}) \quad \textbf{(E.31)}$$

and from the second unnumbered equation preceding **[37, Eq. (1.5)]** we have (see **(11)**)

$$\sum_{l=1}\frac{\psi(\frac{1}{2}l)}{l^3}=-2Li_4(\tfrac{1}{2})-(\gamma+\tfrac{7}{4}\ln(2))\zeta(3)+\frac{1}{120}\pi^4+\tfrac{1}{12}\pi^2\ln^2(2)-\tfrac{1}{12}\ln^4(2). \quad \textbf{(E.32)}$$

By considering its even and odd terms, we also find

$$\sum_{l=1}(-1)^l\frac{\psi(\frac{1}{2}l)}{l^3}=(\tfrac{3}{4}\gamma+\tfrac{7}{4}\ln(2))\zeta(3)-\tfrac{11}{1440}\pi^4-\tfrac{1}{12}\pi^2\ln^2(2)+\tfrac{1}{12}\ln^4(2)+2Li_4(\tfrac{1}{2}).\qquad\textbf{(E.33)}$$

Extracted from **[10**, Section 8.2],

$$\sum_{l=0}(-1)^l\frac{\psi(l)}{l+\frac{1}{2}}=\pi(1-\gamma/2-\ln(2))+2(\gamma+G)-4+2\ln(2)\qquad\textbf{(E.34)}$$

and

$$\sum_{l=0}(-1)^l\frac{\psi(l+1)}{(l+\frac{1}{2})^3}=\tfrac{1}{32}\psi^{(3)}(\tfrac{1}{4})-\frac{\pi^4}{4}-\tfrac{7}{2}\pi\zeta(3)-\frac{\pi^3}{2}\ln(2)-\gamma\frac{\pi^3}{4}.\qquad\textbf{(E.35)}$$

where the trigamma function arises from **[14]** and the imaginary part of *Li₄(i)*. Also see (**C.27**) and **C.28**).

By rewriting, and disentangling the remaining (corrected – also see (**E.26**)) cases presented in **[37**, **Example (3.6)]** (in the fifth case of this example collection, change the coefficient $\frac{41}{16}$ to $\frac{61}{16}$) we find, after considerable effort, the following relations

$$\sum_{l=1}\frac{\psi(\frac{1}{2}l)}{l^2}=-\tfrac{1}{6}\pi^2\gamma-\tfrac{5}{8}\zeta(3),\qquad\textbf{(E.36)}$$

$$\sum_{l=1}\frac{\psi(\frac{1}{2}l)^2}{l^2}=\zeta(3)(\ln(2)-\tfrac{5}{4}\psi(\tfrac{1}{2}))+\tfrac{1}{6}\ln^4(2)+\tfrac{1}{2}\pi^2\ln^2(2)+\tfrac{2}{3}\pi^2\ln(2)\psi(\tfrac{1}{2})$$
$$+\tfrac{1}{6}\pi^2\psi(\tfrac{1}{2})^2-\tfrac{1}{120}\pi^4+4Li_4(\tfrac{1}{2}),\qquad\textbf{(E.37)}$$

$$\sum_{l=1}\frac{\psi(\frac{1}{2}l)\psi(l)}{l^2}=\tfrac{1}{8}\zeta(3)(3\psi(\tfrac{1}{2})-\ln(2))+\tfrac{1}{6}\pi^2\psi(\tfrac{1}{2})^2+\tfrac{2}{3}\pi^2\ln(2)\psi(\tfrac{1}{2})+\tfrac{17}{24}\pi^2\ln^2(2)$$
$$-\tfrac{1}{24}\ln^4(2)+\tfrac{1}{32}\pi^4-Li_4(\tfrac{1}{2}),\qquad\textbf{(E.38)}$$

$$\sum_{l=0}\frac{\psi(\frac{1}{2}+l)^2}{(\frac{1}{2}+l)^2}=\frac{\pi^2}{2}\psi(\tfrac{1}{2})^2+\psi(\tfrac{1}{2})(2\pi^2\ln(2)-7\zeta(3))+\tfrac{4}{3}\pi^2\ln^2(2)$$
$$+\tfrac{2}{3}\ln^4(2)-\tfrac{23}{360}\pi^4+16Li_4(\tfrac{1}{2}),\qquad\textbf{(E.39)}$$

$$\sum_{l=1}\frac{\psi(\frac{1}{2}+l)\psi(l)}{l^2}=\tfrac{9}{2}\zeta(3)\psi(\tfrac{1}{2})+\tfrac{1}{6}\pi^2\psi(\tfrac{1}{2})^2+\tfrac{1}{3}\pi^2\ln(2)\psi(\tfrac{1}{2})+\tfrac{49}{360}\pi^4$$
$$+\tfrac{1}{3}\pi^2\ln^2(2)-\tfrac{1}{3}\ln^4(2)-8Li_4(\tfrac{1}{2}),\qquad\textbf{(E.40)}$$

$$\sum_{l=1}\frac{\psi(\frac{1}{2}+l)\psi(l)}{(l-\frac{1}{2})^2}=\tfrac{21}{2}\zeta(3)\psi(\tfrac{1}{2})+\tfrac{1}{2}\pi^2\psi(\tfrac{1}{2})^2+\pi^2\ln(2)\psi(\tfrac{1}{2})+\tfrac{49}{180}\pi^4$$
$$+\tfrac{2}{3}\pi^2\ln^2(2)-\tfrac{2}{3}\ln^4(2)-16Li_4(\tfrac{1}{2}),\qquad\textbf{(E.41)}$$

and

$$\sum_{l=0}\frac{\psi(l+\frac{1}{2})}{(l+\frac{1}{2})^3}=7\zeta(3)\psi(\tfrac{1}{2})+\tfrac{2}{3}\pi^2\ln^2(2)-\tfrac{2}{3}\ln^4(2)+\tfrac{23}{360}\pi^4-16Li_4(\tfrac{1}{2}).\qquad\textbf{(E.42)}$$

Note that **(E.42)** is required for the proofs of (**E.8**), (**E.31**), (**E.33**), (**E.40**) and (**E.41**). Also, the separation of (**E.38**) into its even and odd components is required for the derivation of (**E.40**) and (**E.41**).

From **[2, Eq. (3.71)]** and **[2, Eq. (3.93)],** among many other sums expressed in terms of harmonic numbers, we find two interesting and representative results

$$\sum_{l=1} \frac{\psi'(l+1)^2}{(l+1)^2} = \tfrac{101}{22680}\pi^6 - \tfrac{1}{36}\pi^4 - \zeta(3)^2 \qquad \textbf{(E.43)}$$

$$\sum_{l=1} \frac{\psi'(l+\frac{1}{2})^2}{(l^2-\frac{1}{4})} = \tfrac{\pi^2}{2}(\pi^2 - 8\zeta(3)) \qquad \textbf{(E.44)}$$

Also, see **(E.24)** and **(E.9).** And finally, from [**38**] (in respective order: Eqs. (3.2a), (2.5c) – incorrect: the negative sign in the summand should be a plus sign – and (2.5d) – incorrect – the right-hand side should be preceded by a minus sign) we have, making use of **(66)** with n=4

$$\sum_{l=1} \frac{\psi(l+1)^2}{l^3} = (\gamma^2 - \tfrac{1}{6}\pi^2)\zeta(3) + \tfrac{7}{2}\zeta(5) - \tfrac{1}{36}\gamma\pi^4 \qquad \textbf{(E.45)}$$

$$\sum_{l=1} \frac{\psi'(l+1)}{l^3} = -\tfrac{1}{3}\pi^2\zeta(3) + \tfrac{9}{2}\zeta(5) \qquad \textbf{(E.46)}$$

and

$$\sum_{l=1} \frac{\psi(l)\psi'(l)}{l^2} = \tfrac{1}{3}\pi^2\zeta(3) - \tfrac{7}{2}\zeta(5) - \tfrac{7}{360}\gamma\pi^4 . \qquad \textbf{(E.47)}$$